\documentclass{amsart}
\newtheorem{thr}{Theorem}[section]
\newtheorem{lem}[thr]{Lemma}

\theoremstyle{definition}
\newtheorem{defn}[thr]{Definition}
\newtheorem{ex}[thr]{Example}
\newtheorem{quest}[thr]{Question}

\theoremstyle{remark}
\newtheorem{remr}[thr]{Remark}

\numberwithin{equation}{section}

\def\trop{rk_t}
\def\Kap{rk_K}
\def\Kf{\textbf{K}}
\def\S{\mathcal S}

\def\B{\mathbb{B}}
\def\C{\mathbb{C}}
\def\R{\mathbb{R}}
\def\i{\infty}
\def\P{\mathcal P}
\def\I{\mathcal I}

\def\T{\mathbb{T}}
\def\ow{\overline{W}}

\def\Su{\textrm{Supp}}

\begin{document}

\title{When do the r-by-r minors of a matrix form a tropical basis?}

\author{Yaroslav Shitov}
\address{Moscow State University, Leninskie Gory, 119991, GSP-1, Moscow, Russia}
\email{yaroslav-shitov@yandex.ru}


\begin{abstract}
We show that the $r$-by-$r$ minors of a $d$-by-$n$ matrix of variables form a tropical basis of the
ideal they generate if and only if $r\leq3$, or $r=\min\{d,n\}$, or else
$r=4$ and $\min\{d,n\}\leq6$. This answers a question asked by M.~Chan, A.~Jensen, and E.~Rubei.
\end{abstract}

\maketitle

\section{Introduction}

The \textit{tropical semiring} is the set $\R$ of real numbers with the operations of tropical addition
and tropical multiplication that are defined by $$a\oplus b=\min\{a,b\},\,\,a\otimes b=a+b.$$
Being important for many different applications (see~\cite{AGG, BCOQ, HOW}),
the tropical semiring is also of considerable interest for studying algebraic geometry, see~\cite{EML, Mikh}.
One of the important concepts is the notion of the rank of a tropical matrix, see~\cite{AGG, DSS}.
In contrast with the case of matrices over a field, there are many different important rank functions
for tropical matrices. Most of these functions have been described and were investigated in~\cite{AGG, DSS}.

It is sometimes useful to extend the tropical semiring with an infinite positive element, denoted by $\i$.
The semiring $\left(\R\cup\{\i\},\min,+\right)$ is called the \textit{completed tropical semiring} and denoted by $\T$.
The binary boolean semiring $\B=\left(\{0,\infty\},\min,+\right)$ will also be useful for our considerations.
We start with the definition of tropical rank, which is one of the most important notions for studying tropical matrices.

\begin{defn}\label{permdef}\label{tropdeg}
The tropical \textit{permanent} of a matrix $S\in\R^{n\times n}$ is defined by
\begin{equation}\label{permanent}\textrm{perm}(S)=\min\limits_{\sigma\in\S_n} \left\{s_{1,\sigma(1)}+\ldots+s_{n,\sigma(n)}\right\},\end{equation}
where $\S_n$ denotes the set of all permutations on $\{1,\ldots,n\}$.
$S$ is called \textit{tropically singular} if the minimum in~(\ref{permanent})
is attained at least twice. Otherwise $S$ is called \textit{tropically non-singular}.
\end{defn}

\begin{defn}\label{troprk}
The \textit{tropical rank}, $\trop(M)$, of a matrix $M$$\in$$\R^{p\times q}$ is the largest number $r$
such that $M$ contains a tropically non-singular $r$-by-$r$ submatrix.
\end{defn}

The notion of tropical linear dependence is also important for our considerations.

\begin{defn}\label{tropdep}
A family of vectors $a_1,\ldots,a_m\in\R^n$ is called \textit{tropically linearly dependent}
(or simply \textit{tropically dependent}) if there exist $\lambda_1,\ldots,\lambda_m\in\T$ such that $(\lambda_1,\ldots,\lambda_m)\neq(\i,\ldots,\i)$
and for every $k\in\{1,\ldots,n\}$ the minimum in $\min_{\tau=1}^m \{\lambda_\tau+a_{\tau k}\}$ is attained at least twice.
In this case, the tuple $(\lambda_1,\ldots,\lambda_m)\in\T^m$ is said to \textit{realize the tropical dependence} of the family $a_1,\ldots,a_m$.
If a family is not tropically dependent, then it is called \textit{tropically independent}.
\end{defn}

The Kapranov rank is the other important notion we deal with.
In order to define this notion we need the following structure, which arises from the study of algebraic geometry.

\begin{defn}\label{Puis}
By $\Kf=\mathbb{C}[[t]]$ we denote the field that consists of the formal sums
$$a(t)=\sum_{e\in\R} a_et^{e}, \mbox{$ $} a_e\in\C,$$ such that the support $E(a)=\{e\in\R: a_e\neq0\}$
is a well-ordered subset of $\R$. The \textit{degree} map $\deg:\Kf^*\rightarrow\R$ takes a sum to the
exponent of its \textit{leading term}, i.e. $\deg a=\min E(a).$ 
By definition, the degree of the zero element of $\Kf$ is $\infty$.
\end{defn}

\begin{remr}
In what follows, the symbol $t$ always denotes the variable of an element of the field $\Kf=\mathbb{C}[[t]]$.
\end{remr}

We can naturally generalize the degree map for matrix arguments. Namely,
we write $B=\deg A$ for matrices $A\in\Kf^{m\times n}$, $B\in\T^{m\times n}$ if $\deg a_{ij}=b_{ij}$ for all $i$, $j$.
In this case, the matrix $A$ is said to be a \textit{lift} of $B$.
The Kapranov rank of a matrix can be defined in the following way (see~\cite[Corollary 3.4]{DSS}).

\begin{defn}\label{Kapr}
The \textit{Kapranov rank} of a matrix $B\in\R^{m\times n}$ is the smallest rank of any lift of $B$, i.e.
$$\Kap(B)=\min\left\{rank(A)\left|\,A\in\Kf^{m\times n}, \deg A=B\right.\right\},$$
where $rank$ is the classical rank function of matrices over the field $\Kf$.
\end{defn}

The notion of Kapranov rank was introduced by Develin, Santos, and Sturmfels in~\cite{DSS}.
They prove the following theorem.

\begin{thr}\label{dssthr}\cite[Theorems 1.4, 5.5, and 6.5]{DSS}
Let $A\in\R^{m\times n}$. Then $\Kap(A)\geq\trop(A)$. If $\trop(A)\leq2$, then $\Kap(A)=\trop(A)$.
If $\trop(A)<\min\{m,n\}$, then $\Kap(A)<\min\{m,n\}$.
\end{thr}

In~\cite{DSS} it was also shown that the functions of tropical and Kapranov rank are indeed different.
The example of a matrix $C\in\R^{7\times7}$ such that $\trop(C)=3$, $\Kap(C)=4$ was provided in~\cite{DSS}.
Also, in~\cite{DSS} the connection between the Kapranov rank and the notion of realizability
of matroids was pointed out.

Further investigations of the Kapranov rank have been carried out in~\cite{CJR, KR, Sh}.
Kim and Roush in~\cite{KR} prove that it is NP-hard to decide whether the Kapranov rank equals $3$ even for $01$-matrices.
It is also shown in~\cite{KR} that there exist matrices with tropical rank $3$ and arbitrarily high Kapranov rank.

Chan, Jensen, and Rubei prove the following theorem.

\begin{thr}\cite[Corollary 1.5]{CJR}\label{cjrthr}
Let a matrix $A\in\R^{d\times n}$ be such that $\min\{d,n\}\leq5$. Then $\trop(A)=\Kap(A)$.
\end{thr}

The following example of a matrix with different tropical and Kapranov ranks is minimal possible
with respect to the size of matrices.

\begin{ex}\label{6x6}\cite[Example 2.1]{Sh}
\textit{Let} $$A=\left(
\begin{array}{cccccc}
0 & 0 & 4 & 4 & 4 & 4 \\
0 & 0 & 2 & 4 & 1 & 4 \\
4 & 4 & 0 & 0 & 4 & 4 \\
2 & 4 & 0 & 0 & 2 & 4 \\
4 & 4 & 4 & 4 & 0 & 0 \\
2 & 4 & 1 & 4 & 0 & 0 \\
\end{array}
\right)
.$$
\textit{Then} $\trop(A)=4$, $\Kap(A)=5$.
\end{ex}

In our paper we use these notions of rank to give an answer for a question
on tropical bases. For a vector $\omega\in\R^n$, the \textit{$\omega$-degree}
of a monomial $x_1^{p_1}\ldots x_n^{p_n}$ is defined to be $\sum_{i=1}^n \omega_ip_i$.
The \textit{initial form} of a polynomial $f\in\C[x_1,\ldots,x_n]$
with respect to $\omega$, $\textrm{in}_\omega(f)$, is the sum of terms in $f$ that have
minimal $\omega$-degree. The tropical hypersurface of $f$ is the set
$${\mathcal T}(f)=\left\{\omega\in\R^n\left|\right.\textrm{in}_\omega(f)\,\mbox{ is not a monomial}\right\}.$$
If $I\subset\C[x_1,\ldots,x_n]$ is an ideal, then the \textit{tropical variety} of $I$ is the set
$${\mathcal T}(I)=\bigcap\limits_{f\in I}{\mathcal T}(f).$$
A \textit{tropical basis} of $I$ is a finite generating set $\{f_1,\ldots,f_m\}$ such that
$${\mathcal T}(I)={\mathcal T}(f_1)\cap\ldots\cap{\mathcal T}(f_m).$$

Let us now consider the $d$-by-$n$ matrix $X$ whose elements are unknowns $x_{ij}$.
By $I^{dn}_r\subset\C[x_{11},\ldots,x_{dn}]$ we denote the ideal that is generated by the
$r$-by-$r$ minors of $X$, by $\mathcal{I}^{dn}_r$ the tropical variety of $I^{dn}_r$.
The fundamental theorem of tropical varieties (see~\cite{SS}) implies
that the Kapranov rank of a matrix $A\in\R^{d\times n}$ is less than $r$ if and only if
$A\in\mathcal{I}^{dn}_r$, see also~\cite{CJR, DSS}.
We also note that the intersection of all tropical hypersurfaces
defined by the $r$-by-$r$ minors of $X$ is exactly the set of all $d$-by-$n$ matrices
whose tropical rank is less than $r$. Thus we can see that the $r$-by-$r$ minors of $X$ form a tropical
basis for $I^{dn}_r$ if and only if every $d$-by-$n$ matrix of tropical rank less than $r$ has the Kapranov
rank less than $r$. The following question was posed in~\cite{CJR}.

\begin{quest}\label{questi}\cite[Question 1.1]{CJR}
For which numbers $d$, $n$, and $r$ do the $(r+1)$-by-$(r+1)$ minors of a $d$-by-$n$ matrix form
a tropical basis? Equivalently, for which $d$, $n$, $r$ does every $d$-by-$n$ matrix of tropical rank at most $r$ have
Kapranov rank at most $r$?
\end{quest}

Our paper gives the answer for Question~\ref{questi}. Namely, we prove the following theorem.

\begin{thr}\label{maintheoremofallthis}
Let $d$, $n$, $r$ be positive integers, $r\leq\min\{d,n\}$. Then the $r$-by-$r$ minors of a $d$-by-$n$ matrix
form a tropical basis if and only if at least one of the following conditions holds:

(1) $r\leq3$;

(2) $r=\min\{d,n\}$;

(3) $r=4$ and $\min\{d,n\}\leq6$.
\end{thr}

\begin{remr}
We assume $r\leq\min\{d,n\}$ in Theorem~\ref{maintheoremofallthis} because
the assumption $r>\min\{d,n\}$ gives a degenerate case of Question~\ref{questi}.
Indeed, by definitions, the empty set always generates the zero ideal and trivially
forms a tropical basis for it.
\end{remr}

The paper is organized as follows. In Section 2, we obtain auxiliary results that are helpful to prove the main results of our paper.
The characterization of the tropical rank via linear dependence, which was proved by Z.~Izhakian in~\cite{Iz2} in rather a complicated
way, is obtained as a corollary of Theorem~\ref{dssthr}. Section 3 is devoted to the proof of the fact that the $4$-by-$4$ minors of a
$6$-by-$n$ matrix form a tropical basis. In Section 4, we finalize the proof of the main result.

\smallskip

The following notation will be used throughout our paper. By $A[r_1,\ldots,r_k]$ we will denote the matrix
which is formed by the rows of a matrix $A$ with numbers $r_1,\ldots,r_k$. We will also use the designation
$A[r_1,\ldots,r_k|c_1,\ldots,c_l]$ for the matrix which is formed by the columns of $A[r_1,\ldots,r_k]$ with
numbers $c_1,\ldots,c_l$. By $A_{ij}$ we will denote the cofactors of an $n$-by-$n$ matrix $A$ over a field,
i.e. $A_{ij}=(-1)^{i+j}\det A[1,\ldots,i-1,i+1,\ldots,n|1,\ldots,j-1,j+1,\ldots,n]$.
Also, we will abbreviate the collocation 'without a loss of generality' by 'w.l.o.g.',
'permutations of rows and columns' by 'p.r.c.'

\section{Preliminary results}

This section provides some auxiliary results that are helpful to prove the main results of our paper.
In particular, we introduce the notion of a pattern of a matrix and obtain straightforward properties of definitions.
We also consider systems of a small number of linear equations over $\Kf$.

\subsection{The pattern and simple properties of tropical matrices}

First, we introduce the notion of a pattern of a matrix.

\begin{defn}\label{patt}
The \textit{pattern} of a matrix $A\in\R^{m\times n}$ is the matrix $B\in\B^{m\times n}$
whose entries $(b_{uv})$ are defined by
$$b_{uv}=\begin{cases}
0 \mbox{ if } a_{uv}=\min_{i=1}^m\{a_{iv}\},\\
\i \mbox{ if } a_{uv}>\min_{i=1}^m\{a_{iv}\}.
\end{cases}$$
The pattern of $A$ is denoted by $\P(A)$.
\end{defn}

\begin{defn}\label{support}
The \textit{support} of a vector $b\in\B^{n}$ is the set of all $i\in\{1,\ldots,n\}$ such that
the $i$th coordinate of $b$ equals $0$, i.e. $b_{i}=0$.
The support of $j$th column of a matrix $A\in\B^{m\times n}$ is denoted by $\Su_j(A)$.
\end{defn}

The concept of tropical rank is also useful in the case of matrices over $\B$.

\begin{defn}\label{Btrop1}
A family of vectors $a_1,\ldots,a_m\in\B^n$ is called \textit{$\B$-tropically linearly dependent}
(or simply \textit{$\B$-tropically dependent}) if there exists a nonempty set $I\subset\{1,\ldots,m\}$
such that for every $k\in\{1,\ldots,n\}$ the cardinality of the set $\{i\in I\left|\right.a_{ik}=0\}$ is different from $1$.
In this case, $I$ is said to \textit{realize the $\B$-tropical dependence} of the family $a_1,\ldots,a_m$.
If a family is not $\B$-tropically dependent, then it is called \textit{$\B$-tropically independent}.
\end{defn}

\begin{lem}\label{r->b}
Let $A\in\R^{m\times n}$, $B=\P(A)$. Let numbers $\{\lambda_\tau\}$ ($\tau$ runs over a set $\I\subset\{1,\ldots,m\}$)
realize the tropical dependence of the rows of $A$ with numbers from $\I$. Set
$I=\{\rho\in\I|\lambda_\rho=\min_{\tau\in\I}\{\lambda_\tau\}\}$.
Then $I$ realizes the $\B$-tropical dependence of the rows of $B$.
\end{lem}

\begin{proof}
Assume the converse. Then, by Definition~\ref{Btrop1}, for some
$k\in\{1,\ldots,n\}$, $i\in I$ it holds that $a_{ik}=0$, and for every $i'\in I\setminus\{i\}$ it holds that $a_{i'k}=\i$.
In this case, by Definition~\ref{patt}, $\lambda_j+a_{jk}>\lambda_i+a_{ik}$
for every $j\in\I\setminus\{i\}$. This contradicts Definition~\ref{tropdep}.
\end{proof}

\begin{defn}\label{Btrop3}
Let $n>1$. A matrix $S\in\B^{n\times n}$ is called \textit{$\B$-tropically singular} if the minimum in~(\ref{permanent})
is attained at least twice. Otherwise $S$ is called \textit{$\B$-tropically non-singular}.
\end{defn}

\begin{lem}\label{Btrop}
The size of any $\B$-tropically non-singular submatrix of the pattern of a matrix $A\in\R^{n\times m}$
cannot exceed the tropical rank of $A$.
\end{lem}

\begin{proof}
By Definition~\ref{patt}, $\B$-tropical non-singularity of a submatrix $\P(A)[r_1,\ldots,r_k|c_1,\ldots,c_k]$
implies the tropical non-singularity of the submatrix $A[r_1,\ldots,r_k|c_1,\ldots,c_k]$.
Now the lemma follows from Definition~\ref{troprk}.
\end{proof}

We introduce the natural equivalence relation for tropical matrices.

\begin{defn}\label{ekviv}
The \textit{equivalent transformations} of a matrix $A\in\R^{m\times n}$ are the permutations
of rows, of columns and the tropical multiplication of some row or column by a number $u\in\R$.
Matrices $A,B\in\R^{m\times n}$ are said to be \textit{equivalent} if they can be obtained from each other
by applying a composition of equivalent transformations.
\end{defn}

The following lemmas follow directly from definitions.

\begin{lem}\label{multiply}
If matrices $A,B\in\R^{m\times n}$ are equivalent, then $\trop(A)=\trop(B)$, $\Kap(A)=\Kap(B)$.
\end{lem}

\begin{lem}\label{monot}
Let $A\in\R^{m\times n}$. If $B$ is a submatrix of $A$, then $\trop(A)\geq\trop(B)$, $\Kap(A)\geq\Kap(B)$.
\end{lem}

\begin{lem}\label{conve}
Let tuples $(\lambda_1,\ldots,\lambda_m)$, $(\mu_1,\ldots,\mu_m)$ both realize the tropical dependence of a family $a_1,\ldots,a_m\in\R^n$.
Then the tuple $(\min\{\lambda_1,\mu_1\},\ldots,\min\{\lambda_m,\mu_m\})$ realizes the tropical dependence of the family $a_1,\ldots,a_m$ as well.
\end{lem}

\begin{lem}\label{li->tr}
Let $A\in(\Kf^*)^{m\times n}$. Let a nonzero vector $(\lambda_1,\ldots,\lambda_m)\in\Kf^m$ be such that $\sum_{i=1}^m\lambda_{i} a_{ik}=0$
for every $k\in\{1,\ldots,n\}$. Then the tuple
$(\deg\lambda_1,\ldots,\deg\lambda_m)\in\T^m$ realizes the tropical dependence of the matrix $\deg A\in\R^{m\times n}$.
\end{lem}

\begin{proof}
Assume the converse. Then for some $k\in\{1,\ldots,n\}$ it holds that $\deg\lambda_{i_0}+\deg a_{i_0k}<\deg\lambda_{i}+\deg a_{ik}$
for every $i\in\{1,\ldots,m\}\setminus\{i_0\}$.
This implies that $\deg\left(\sum_{i=1}^m\lambda_i a_{ik}\right)=\deg\lambda_{i_0}+\deg a_{i_0k}\neq\i$ and gives a contradiction.
\end{proof}

\subsection{Systems of linear equations over $\Kf$}

It will be very important for our further considerations to decide whether a given system of linear equations over $\Kf$
has a solution with prescribed degrees of unknowns. This section provides some important sufficient conditions
for systems with three, two, or a single equation.

\begin{lem}\label{sys}
Let $A\in\Kf^{6\times3}$, $(h_1,\ldots,h_6)\in\R^6$, $$D=\sum_{i=1}^3\min_{j=1}^6\{\deg a_{ji}+h_j\}<\i.$$ Let
numbers $u,v,y,z\in\{1,2,3,4,5,6\}$ be pairwise distinct. Let $\deg\det A[p,q,r]=D-(h_p+h_q+h_r)$
for every pairwise distinct $p,q,r\in\{u,v,y,z\}$. Then there exist $x_1,\ldots,x_6\in\Kf$ such that
$\deg x_j=h_j$ for every $j\in\{1,2,3,4,5,6\}$, and $\sum_{j=1}^6 a_{ji}x_j=0$ for every $i\in\{1,2,3\}$.
\end{lem}

\begin{proof}
1. We assume w.l.o.g. that $\{u,v,y,z\}=\{1,2,3,4\}$. We set $x_4=\xi t^{h_4}$, $x_5=t^{h_5}$,
$x_6=t^{h_6}$, where $\xi$ is a certain element of $\C^*$.

2. Denote
\begin{equation}\label{prel1}A'=\left(
\begin{array}{ccc}
a_{11}&a_{12}&a_{13}\\
a_{21}&a_{22}&a_{23}\\
a_{31}&a_{32}&a_{33}\\
\sum_{q=4}^6 a_{q1}x_q&\sum_{q=4}^6 a_{q2}x_q&\sum_{q=4}^6 a_{q3}x_q\\
\end{array}
\right)
\end{equation}
We define the elements $x_1$, $x_2$, $x_3$ to be a (unique) solution of the linear system
$\sum_{p=1}^3 a'_{pi}x_p=a'_{4i},$ $i\in\{1,2,3\}$.
It remains to check that $\deg x_j=h_j$ for any $j\in\{1,2,3,4,5,6\}$.

3. From item 1 it follows directly that $\deg x_{j'}=h_{j'}$ for any $j'\in\{4,5,6\}$.

4. Further, assume 
$j=1$. The Cramer's rule for solving linear systems shows that $$x_1=\frac{\det A'[4,2,3]}{\det A'[1,2,3]}.$$
By the assumptions of the lemma, $\deg\det A[1,2,3]=D-(h_1+h_2+h_3)$. Thus we need to check
$\deg\det A'[2,3,4]=D-(h_2+h_3)$. By~(\ref{prel1}),
$$\det A'[2,3,4]=x_4\det A[2,3,4]+(a_{51}x_5+a_{61}x_6)(a_{22}a_{33}-a_{23}a_{32})+$$
\begin{equation}\label{prel2}+(a_{52}x_5+a_{62}x_6)(a_{23}a_{31}-a_{21}a_{33})+(a_{53}x_5+a_{63}x_6)(a_{21}a_{32}-a_{22}a_{31}).\end{equation}

Let us now obtain a lower bound for the degrees of the terms in the right-hand side of~(\ref{prel2}).
For any $\alpha\in\{4,5,6\}$, $\beta,\gamma\in\{1,\ldots,6\}$ it holds that
$$\deg (x_\alpha a_{\alpha1}a_{\beta2}a_{\gamma3})=\deg\left(x_\alpha a_{\alpha1}\right)+\deg a_{\beta2}+\deg a_{\gamma3}\geq$$
$$\geq\min_{j=1}^6\{a_{j1}+h_j\}+\left(\min_{j=1}^6\{a_{j2}+h_j\}-h_2\right)+\left(\min_{j=1}^6\{a_{j3}+h_j\}-h_3\right)=D-(h_2+h_3).$$
By item 1, $\deg$$x_4$$=$$h_4$, so the assumption of the lemma implies that $\deg$$\left(x_4\det A[2,3,4]\right)$$=$$D$$-$$($$h_2$$+$$h_3$$)$.
This shows that the term $x_4\det A[2,3,4]=\xi t^{h_4}\det A[2,3,4]$ has the lowest degree among the
terms of the right-hand side of~(\ref{prel2}). Therefore the condition $\deg\det A'[2,3,4]=D-(h_2+h_3)$,
and thus the condition $\deg x_1=h_1$, holds for all but one complexes $\xi$.

5. In the same way we can prove that the conditions $\deg x_2=h_2$ and
$\deg x_3=h_3$ also hold for all but one or two complex numbers $\xi$.
\end{proof}

We need the following lemma to prove a similar statement for systems with two equations.

\begin{lem}\label{forsystem2}
Let $S\in\Kf^{2\times m}$. Then there exists $\xi\in\mathbb{C}^*$ such that
$\deg(\xi s_{1k}+s_{2k})=\min\{\deg s_{1k}, \deg s_{2k}\}$ for every $k\in\{1,\ldots,m\}$.
\end{lem}

\begin{proof}
If $s_{ik}\neq0$, we denote the coefficient of the leading term of $s_{ik}$ by $\sigma_{ik}$.
If $s_{ij}=0$, we choose $\sigma_{ik}\in\mathbb{C}^*$ arbitrarily. Now it remains to choose
$\xi\in\mathbb{C}\setminus\{0,-\frac{\sigma_{21}}{\sigma_{11}},\ldots,-\frac{\sigma_{2m}}{\sigma_{1m}}\}$.
\end{proof}

\begin{lem}\label{sys2}
Let each column of a matrix $A\in\Kf^{5\times2}$ contain a nonzero element,
and $\deg(a_{p1}a_{q2}-a_{q1}a_{p2})=\min\{\deg a_{p1}+\deg a_{q2}, \deg a_{q1}+\deg a_{p2}\}$ for any different $p,q\in\{1,\ldots,5\}$.
Let $(h_1,\ldots,h_5)\in\R^5$, we denote by $\Theta_i$ ($i\in\{1,2\}$) the set of all $\eta\in\{1,2,3,4,5\}$
that deliver the minimum for $\min_\eta\{\deg a_{\eta i}+h_\eta\}$.
Let also $|\Theta_1|\geq2$, $|\Theta_2|\geq2$, $|\Theta_1\cup\Theta_2|\geq3$.
Then there exist $x_1,\ldots,x_5\in\Kf$ such that $\deg x_j=h_j$ for every $j\in\{1,2,3,4,5\}$,
and $\sum_{j=1}^5 a_{ji}x_j=0$ for $i\in\{1,2\}$.
\end{lem}

\begin{proof}
We denote $\theta_1=\min_{j=1}^5\{a_{j1}+h_j\}$, $\theta_2=\min_{j=1}^5\{a_{j2}+h_j\}$.
We assume w.l.o.g. that $1\in\Theta_1$, $2\in\Theta_2$, and both $\Theta_1$ and $\Theta_2$ have a non-empty
intersection with $\{3,4,5\}$. These settings imply that
$\min_{\iota=3}^5\left\{\deg\det A[1,\iota]+h_1+h_{\iota}\right\}=\min_{\iota=3}^5\left\{\deg\det A[2,\iota]+h_2+h_{\iota}\right\}=\theta_1+\theta_2$.
From Lemma~\ref{forsystem2} it then follows that there exist $\xi_3,\xi_4,\xi_5\in\mathbb{C}^*$ such that
$$
\deg\left(\sum_{\iota=3}^5\det A[1,\iota]t^{h_{1}+h_{\iota}}\xi_\iota\right)=
\deg\left(\sum_{\iota=3}^5\det A[2,\iota]t^{h_{2}+h_{\iota}}\xi_\iota\right)=\theta_1+\theta_2.
$$
From Cramer's rule it then follows that the solution $(y_1,y_2)$ of
\begin{equation}\label{eqsys2-2}\begin{cases}
a_{11}t^{h_{1}}y_1+a_{21}t^{h_{2}}y_2=-\sum_{\iota=3}^5\xi_\iota a_{\iota1}t^{h_{\iota}},\\
a_{12}t^{h_{1}}y_1+a_{22}t^{h_{2}}y_2=-\sum_{\iota=3}^5\xi_\iota a_{\iota2}t^{h_{\iota}}
\end{cases}\end{equation}
is such that $\deg y_1=\deg y_2=0$. We set $x_{1}=y_1t^{h_{1}}$, $x_{2}=y_2t^{h_{2}}$,
$x_{\iota}=\xi_\iota t^{h_{\iota}}$ for $\iota\in\{3,4,5\}$. The equations~(\ref{eqsys2-2})
imply that $\sum_{j=1}^5 a_{j1}x_{j}=\sum_{j=1}^5 a_{j2}x_{j}=0$.
\end{proof}

We also prove a similar lemma for a single linear equation.

\begin{lem}\label{sys1e}
Let $l\in\Kf^{m}$ be such that $\deg l$ realizes the tropical dependence of a vector $a\in\R^m$.
Then there exist $x_1,\ldots,x_m\in\Kf$ such that $\deg x_j=a_j$ for every $j\in\{1,\ldots,m\}$,
and $\sum_{j=1}^m x_jl_j=0$.
\end{lem}

\begin{proof}
Let the minimum in $\min_{j=1}^m \{a_j+\deg l_j\}$ be provided by $j_1,\ldots,j_k\in\{1,\ldots,m\}$.
Definition~\ref{tropdep} implies that $k>1$. Now it is enough to set $x_\jmath=t^{a_\jmath}$
for $\jmath\in\{1,\ldots,m\}\setminus\{j_1\}$, and $x_{j_1}=-\frac{\sum_{\jmath}l_\jmath x_\jmath}{l_{j_1}}$.
\end{proof}

The following equivalent characterization of the tropical rank was proved by Izhakian in~\cite{Iz2}
in rather a complicated way, see also~\cite{AGG}. We now obtain this characterization as a corollary of Theorem~\ref{dssthr}.

\begin{thr}\label{tropthr}
The tropical rank of a matrix $A\in\R^{m\times n}$ equals the cardinality of the largest tropically independent family of rows of $A$.
\end{thr}

\begin{proof}
Let $r_1,\ldots,r_c$ be the numbers of the rows the largest tropically independent family of $A$. Lemma~\ref{sys1e} now shows that
the Kapranov rank of any submatrix $A[r'_1,\ldots,r'_{c+1}]$ is at most $c$, and Theorem~\ref{dssthr} thus shows that
$\trop(A)\leq c$.

On the other hand, if $\trop(A)<c$, then Theorem~\ref{dssthr} implies that $\Kap(A[r_1,\ldots,r_c])<c$. In this case,
Theorem~\ref{li->tr} shows that the rows of $A[r_1,\ldots,r_c]$ are tropically dependent, so indeed $\trop(A)=c$.
\end{proof}

\section{The $4$-by-$4$ minors of a $6$-by-$n$ matrix form a tropical basis}

In this section, we prove that every $6$-by-$n$ matrix of tropical rank $3$ has Kapranov rank $3$ as well.
First, we give a characterization of $6$-by-$n$ matrices with tropical rank $3$ and greater Kapranov rank.
This characterization reduces the problem to a number of special cases each of which can be treated separately.

\subsection{Reduction of the problem to a number of special cases}

The following theorem gives a characterization of $6$-by-$n$ matrices with tropical rank $3$ and greater Kapranov rank.

\begin{thr}\label{hugecases0}
Let a matrix $A\in\R^{6\times n}$ be such that $\trop(A)=3$, $\Kap(A)>3$. Then the equivalence class of $A$ contains a matrix $W$
such that every row of $\P(W)$ contains at least one $0$, and at least one of the following conditions holds:

(i) $\P(W)$ consists of the columns with the supports $\{1,2\}$, $\{3,4\}$, $\{5,6\}$;

(ii) $\P(W)$ contains columns with the support $\{1,2\}$, contains at least one of columns
with the supports $\{4,5,6\}$, $\{3,4,6\}$, $\{3,5,6\}$, $\{3,4,5\}$. $\P(W)$ may also contain
the column with the support $\{3,4,5,6\}$;

(iii) $\P(W)$ may contain several columns with the supports $\{1,2,3\}$, $\{1,4,5\}$, $\{2,4,6\}$, $\{3,5,6\}$.
The supports of other columns of $\P(W)$ have cardinality at least $4$.
If the support of $\jmath$th column of $\P(W)$ has cardinality $4$,
then no proper subset of $\Su_\jmath(\P(W))$ is a support of a column of $\P(W)$;

(iv) $\P(W)$ contains columns with the supports $\{1,2,3\}$, $\{4,5,6\}$,
the support of any other column of $\P(W)$ contains at least two numbers from $\{1,2,3\}$ and at least two numbers from $\{4,5,6\}$;

(v) The set of the supports of the columns of $\P(W)$ is $\left\{\{1,2\},\{1,2,3\},\{4,5,6\}\right\}$.
\end{thr}

\begin{proof}
1. Let $A\in\R^{6\times n}$ be such that $\trop(A)=3$, $\Kap(A)>3$.
Then by Theorem~\ref{tropthr}, the rows of $A$ are tropically dependent.
We apply Definition~\ref{tropdep}. So by Lemma~\ref{multiply}, we assume w.l.o.g. that every column of $\P(A)$
contains at least two zeros.

2. Moreover, by Lemma~\ref{multiply}, we assume w.l.o.g. that any row of $\P(A)$
contains at least one zero. The situation splits into the following three cases.

\textbf{Case O.} Let every column of $\P(A)$ contain at least $4$ zeros.
This case satisfies the assumptions of item (iii) of the theorem we prove.

\textbf{Case A.} Let every column of $\P(A)$ contain at least $3$ zeros, and
some column of $\P(A)$ contain exactly $3$ zeros. Let us consider $2$ special cases, A1 and A2.

A1. Assume there are $u',v'\in\{1,\ldots,n\}$ such that $\left|\Su_{u'}(\P(A))\right|=3$,
$\left|\Su_{v'}(\P(A))\setminus\Su_{u'}(\P(A))\right|=1$. In this case, up to p.r.c.
$$\P(A)=\left(\begin{array}{cc|ccc}
0&a_{12}&&&\\
0&a_{22}&&\P(A)[1,2,3,4|3,\ldots,n]&\\
0&a_{32}&&&\\
\i&0&&&\\\hline
\i&\i&a_{53}&\ldots&a_{5n}\\
\i&\i&a_{63}&\ldots&a_{6n}\end{array}\right).$$

If the $5$th and $6$th columns of $\P(A)$ are different, then by Lemma~\ref{Btrop}, $\trop(A)\geq4$.
This contradiction shows that the $5$th and $6$th columns of $\P(A)$ coincide.
By item 2 any row of $\P(A)$ contains a zero, so up to p.r.c. we obtain
$$\P(A)=\left(\begin{array}{cc|ccc|ccc}
0&p_{12}&&&&\\
0&p_{22}&&\overline{P}&&&\overline{P'}&\\
0&p_{32}&&&&\\
\i&0&&&&\\\hline
\i&\i&\i&\ldots&\i&0&\ldots&0\\
\i&\i&\i&\ldots&\i&0&\ldots&0\end{array}\right).$$
The assumption of Case A shows that any column of the matrix $\overline{P}$ contains at least $3$ zeros.

Now we add a small enough (with respect to the absolute value) $-\varepsilon<0$ to every element of the $5$th and $6$th rows of $A$.
Definition~\ref{patt} shows that the pattern of the matrix obtained equals
$$\left(\begin{array}{cc|ccc|ccc}
0&p_{12}&&&&\i&\ldots&\i\\
0&p_{22}&&\overline{P}&&\ldots&&\ldots\\
0&p_{32}&&&&\ldots&&\ldots\\
\i&0&&&&\i&\ldots&\i\\\hline
\i&\i&\i&\ldots&\i&0&\ldots&0\\
\i&\i&\i&\ldots&\i&0&\ldots&0\end{array}\right)$$
and up to p.r.c. satisfies the assumptions of case (ii).

A2. Now assume that for some $u'',v''\in\{1,\ldots,n\}$ the supports $\Su_{u''}(\P(A))$ and $\Su_{v''}(\P(A))$ are disjoint.
In this case $A$ up to p.r.c. satisfies the assumptions of either item A1 or case (iv).

A3. Let us conclude the analysis of Case A.
Item A1 shows that if some support with cardinality $4$ includes $\Su_u(\P(A))$ for some $u\in\{1,\ldots,n\}$,
then the statement of the theorem holds for $A$. So we can assume w.l.o.g. that no proper subset of a support of cardinality $4$
is a support of a column of $\P(W)$.

By items A1 and A2, the statement of the theorem also holds for $A$ if $\left|\Su_{\widehat{u}}(\P(A))\cap\Su_{\widehat{v}}(\P(A))\right|\in\{0,2\}$
for some $\widehat{u},\widehat{v}\in\{1,\ldots,n\}$ such that $\left|\Su_{\widehat{u}}(\P(A))\right|=\left|\Su_{\widehat{v}}(\P(A))\right|=3$.
Thus we can also assume w.l.o.g. that $\left|\Su_{\widehat{u}}(\P(A))\cap\Su_{\widehat{v}}(\P(A))\right|\in\{1,3\}$ for every such
$\widehat{u}$, $\widehat{v}$. Now it is straightforward to see that up to p.r.c. $A$ satisfies the assumptions of item (iii).

\textbf{Case B.} Let some column of $\P(A)$ contain exactly $2$ zeros.
Assume w.l.o.g.
$$\left(\begin{array}{c}
0\\
0\\
\i\\
\i\\
\i\\
\i\end{array}\right)$$
is a column of $\P(A)$.
The situation splits into the following $3$ cases.

B1. Assume that some column of the matrix $\P(A)[3,4,5,6]$ contains exactly $2$ zeros.
Then by Lemma~\ref{Btrop}, up to p.r.c.
$$\P(A)=\left(\begin{array}{ccc|ccc|ccc}
0&\ldots&0&&P'_1&&&P'_2&\\
0&\ldots&0&&&&&&\\\hline
\i&\ldots&\i&0&\ldots&0&&P'_3&\\
\i&\ldots&\i&0&\ldots&0&&&\\\hline
\i&\ldots&\i&\i&\ldots&\i&0&\ldots&0\\
\i&\ldots&\i&\i&\ldots&\i&0&\ldots&0\end{array}\right).$$
Now we add a small enough $\varepsilon>0$ to every element of the $3$rd and $4$th rows of $A$,
add $2\varepsilon$ to every element of the $1$st and $2$nd rows of $A$.
Definition~\ref{patt} shows that the pattern of the matrix obtained equals
$$\left(\begin{array}{ccc|ccc|ccc}
0&\ldots&0&\i&\ldots&\i&\i&\ldots&\i\\
0&\ldots&0&\i&\ldots&\i&\i&\ldots&\i\\\hline
\i&\ldots&\i&0&\ldots&0&\i&\ldots&\i\\
\i&\ldots&\i&0&\ldots&0&\i&\ldots&\i\\\hline
\i&\ldots&\i&\i&\ldots&\i&0&\ldots&0\\
\i&\ldots&\i&\i&\ldots&\i&0&\ldots&0\end{array}\right),$$
and up to p.r.c. satisfies the assumptions of case (i).

B2. Now assume that every column of $\P(A)[3,4,5,6]$ contains either no or at least three zeros.
In this case
$$\P(A)\left(\begin{array}{ccc|ccc}
0&\ldots&0&&P'&\\
0&\ldots&0&&&\\\hline
\i&\ldots&\i&&&\\
\i&\ldots&\i&&P''&\\
\i&\ldots&\i&&&\\
\i&\ldots&\i&&&\end{array}\right),$$
where any column of $P''$ contains at least three zeros.

Now we add a small enough $\varepsilon>0$ to every element of the $1$st and $2$nd rows of $A$.
Definition~\ref{patt} shows that the pattern of the matrix $A_3$ obtained is
$$\left(\begin{array}{ccc|ccc}
0&\ldots&0&\i&\ldots&\i\\
0&\ldots&0&\i&\ldots&\i\\\hline
\i&\ldots&\i&&&\\
\i&\ldots&\i&&P''&\\
\i&\ldots&\i&&&\\
\i&\ldots&\i&&&\end{array}\right).$$
We consider the two cases.

B2.1. If $\i$ appears as an entry of $P''$, then $\P(A_3)$ satisfies the assumptions of item (ii).

B2.2. Now we assume that
\begin{equation}\label{equamatr2}\P(A_3)=\left(\begin{array}{c|c}
0\ldots0&\i\ldots\i\\
0\ldots0&\i\ldots\i\\\hline
\i\ldots\i&0\ldots0\\
\i\ldots\i&0\ldots0\\
\i\ldots\i&0\ldots0\\
\underbrace{\i\ldots\i}_{p'}&\underbrace{0\ldots0}_{q'}\end{array}\right),\mbox{ where $p',q'>0$}.\end{equation}
We can assume w.l.o.g. that $p'$ is minimal over all matrices $A''$ that are equivalent to $A$ and such that
the pattern of $A''$ has the form~(\ref{equamatr2}). By Lemma~\ref{multiply}, we can assume w.l.o.g.
that the minimal element of any row of $A_3$ equals $0$.

We add the minimal element of the matrix $\P(A)[3,4,5,6|1,\ldots,p']$
to every element of the $1$st and $2$nd rows of $A_3$.
Denote by $A_3'$ the matrix obtained. By Definition~\ref{patt}, up to p.r.c.
$$\P(A'_3)=\left(\begin{array}{c|c|c}
0\ldots0&0\ldots0&\i\ldots\i\\
0\ldots0&0\ldots0&\i\ldots\i\\\hline
\i\ldots\i&&0\ldots0\\
\i\ldots\i&P_3&0\ldots0\\
\i\ldots\i&&0\ldots0\\
\underbrace{\i\ldots\i}_{p''}&\underbrace{\,\,\,\,\,\,\,\,\,\,\,\,\,\,\,\,\,\,\,\,\,}_{p'-p''}&0\ldots0\end{array}\right),$$
where $p''<p'$, every column of $P_3$ contains at least one zero.

If $p''=0$, then $A_3'$ satisfies the assumptions of either Case O or Case A.
Further we assume that $p''>0$.
The minimality of $p'$ shows that $\i$ appears as an entry of $P_3$.
If every column of $A'_3$ contains at least $2$ zeros, then $A_3'$ satisfies the assumptions of either item B1 or item B2.1.

Thus it remains to consider the case when some column of $P_3$ contains exactly $1$ zero.
By Lemma~\ref{Btrop}, up to p.r.c.
$$\P(A_3')=\left(\begin{array}{c|c|c}
0\ldots0&0\ldots0&P_8\\
0\ldots0&0\ldots0&\\\hline
\i\ldots\i&0\ldots0&0\ldots0\\
\i\ldots\i&\i\ldots\i&0\ldots0\\
\i\ldots\i&\i\ldots\i&0\ldots0\\
\i\ldots\i&\i\ldots\i&0\ldots0\end{array}\right).$$
Now we add a small enough $\varepsilon>0$ to every element of the first three rows of $A_3'$.
Definition~\ref{patt} shows that the matrix obtained to satisfies the conditions of item (v).
This completes the consideration of item B2.2. Note that items B2.1 and B2.2
cover all possible cases, thus the consideration of item B2 is also complete.

B3. Finally, assume that some column of $\P(A)[3,4,5,6]$ contains exactly $1$ zero.
Lemma~\ref{Btrop} shows that up to p.r.c. it holds that
\begin{equation}\label{equamatr1}
\P(A)=\left(\begin{array}{c|c|c}
0\ldots0&P_0&P'_0\\
0\ldots0&&\\\hline
\i\ldots\i&0\ldots0&P''_0\\\hline
\i\ldots\i&\i\ldots\i&0\ldots0\\
\i\ldots\i&\i\ldots\i&0\ldots0\\
\underbrace{\i\ldots\i}_{p}&\underbrace{\i\ldots\i}_{q}&\underbrace{0\ldots0}_r\end{array}\right),\mbox{ $ $} p>0,\mbox{ $ $} q>0,\mbox{ $ $} r>0,
\end{equation}
where any column of $P_0$ contains at least $1$ zero.

\textit{We choose a matrix (and denote it by $A_0$) for which the value $p+q$ is minimal under the following assumptions:
the pattern of $A_0$ has the form~(\ref{equamatr1}), $A_0$ and $A$ are equivalent, every column of $\P(A_0)$ contains at least $2$ zeros.}


By Lemma~\ref{multiply}, we can assume w.l.o.g. that the minimal element of any column of $A_0$ equals $0$.
The matrix $A_0$ is nonnegative and equals
$$\left(\begin{array}{ccc|ccc|ccc}
0&\ldots&0&&A_0[1,2|p+1,\ldots,p+q]&&&A_1&\\
0&\ldots&0&&&&&&\\\hline
&A_0[3|1,\ldots,p]&&0&\ldots&0&&A_2&\\\hline
&&&&&&0&\ldots&0\\
&A_0[4,5,6|1,\ldots,p]&&&A_0[4,5,6|p+1,\ldots,p+q]&&0&\ldots&0\\
&&&&&&0&\ldots&0\end{array}\right).$$
Definition~\ref{patt} implies that the entries of the matrices
$A_0[3|1,\ldots,p]$ and $A_0[4,5,6|1,\ldots,p+q]$ are all positive, 
every column of $A_0[1,2|p+1,\ldots,p+q]$ contains at least $1$ zero.
In particular, the minimal element of $A_0[4,5,6|1,\ldots,p+q]$ (we denote it by $m$) is positive.

We add $m$ to every element of the first three rows of $A_0$, add $-m$
to every element of the first $p+q$ columns of $A_0$. We obtain the nonnegative matrix $C$ that is equal to
$$\left(\begin{array}{ccc|ccc|ccc}
0&\ldots&0&&A_0[1,2|p+1,\ldots,p+q]&&&C_1&\\
0&\ldots&0&&&&&&\\\hline
&A_0[3|1,\ldots,p]&&0&\ldots&0&&C_2&\\\hline
&&&&&&0&\ldots&0\\
&C[4,5,6|1,\ldots,p]&&&C[4,5,6|p+1,\ldots,p+q]&&0&\ldots&0\\
&&&&&&0&\ldots&0\end{array}\right),$$
and $0$ appears as an entry of the matrix $C[4,5,6|1,\ldots,p+q]$,
the entries of the matrices $C_1$ and $C_2$ are all positive.

If every column of $C$ contains at least $3$ zeros, then $C$ satisfies the assumptions of either Case O or Case A.
Thus we can assume that some column of $C$ contains exactly $2$ zeros.
Thus some permutations of the first three rows and of the first $p+q$ columns of $C$ produce the matrix
\begin{equation}\label{equamatr22}C'=\left(\begin{array}{c|c|ccc|ccc}
0&c'_{12}&&&&&&\\
0&c'_{22}&&C'[1,2,3|3,\ldots,p+q]&&&C'_0&\\
\omega_3&0&&&&&&\\\hline
\omega_4&c'_{42}&&&&0&\ldots&0\\
\omega_5&c'_{52}&&C'[4,5,6|3,\ldots,p+q]&&0&\ldots&0\\
\omega_6&c'_{62}&&&&0&\ldots&0\end{array}\right),\end{equation}
where $\omega_\iota>0$ for $\iota\in\{3,4,5,6\}$, every column of $C'[1,2,3|2,\ldots,p+q]$
contains at least $2$ zeros, the minimal element of $C'[4,5,6|2,\ldots,p+q]$ equals $0$,
the elements of $C'_0$ are positive. The situation splits into the following three cases.

B3.1. If every column of the matrix $C'[3,4,5,6|2,\ldots,p+q]$ contains either no or at least $3$ zeros, then the matrix
$C'$ satisfies the assumptions of item B2.

B3.2. If some column of $C'[3,4,5,6|2,\ldots,p+q]$ contains exactly $2$ zeros,
then $C'$ satisfies the assumptions of item B1.

B3.3. Finally, let some column (we denote its number by $j$) of $C'[3,4,5,6|2,\ldots,p+q]$ contain exactly $1$ zero,
i.e. $c'_{ij}=0$, $i\geq3$.

Assume $i\neq3$. Fix an arbitrary  $i_0\in\{4,5,6\}\setminus\{i\}$. Then $c'_{i_0j}>0$, $c'_{3j}>0$.
By~(\ref{equamatr22}), $c'_{3,p+q+1}>0$. Thus the matrix
$$C'[2,i,i_0,3|1,j,p+q+1,2]=\left(\begin{array}{cccc}
0&x_1&x_2&x_3\\
\omega_i&0&x_4&x_5\\
\omega_{i_0}&c'_{i_0j}&0&x_6\\
\omega_{3}&c'_{3j}&c'_{3,p+q+1}&0\\
\end{array}\right),
$$ where $x_1,\ldots,x_6$ are nonnegative, is tropically non-singular.
This implies $\trop(A)>3$ and contradicts the assumptions of the theorem.

Thus $i=3$. We permute the $2$nd and $j$th columns of $C'$ to obtain the matrix $D$ such that
\begin{equation}\label{equamatr3}\P(D)=\left(\begin{array}{c|c|ccc|ccc}
0&c'_{1j}&&&&&&\\
0&c'_{2j}&&\P(D)[1,2,3|3,\ldots,p+q]&&&C'_0&\\
\i&0&&&&&&\\\hline
\i&\i&&&&0&\ldots&0\\
\i&\i&&\P(D)[4,5,6|3,\ldots,p+q]&&0&\ldots&0\\
\i&\i&&&&0&\ldots&0\end{array}\right).\end{equation}
By Lemma~\ref{Btrop},
the rows of $\P(D)[4,5,6]$ coincide. Note that the matrix $D[4,5,6|2,\ldots,p+q]$
equals $C'[4,5,6|2,\ldots,p+q]$ up to the permutation of columns, so
the minimal element of $D[4,5,6|3,\ldots,p+q]$ equals $0$.
Thus the matrix $\P(D)[4,5,6|3,\ldots,p+q]$ contains several columns consisting of zeros,
all the other columns of $\P(D)[4,5,6|3,\ldots,p+q]$ consist of $\i$.
By~(\ref{equamatr3}), this contradicts the minimality of $p+q$.
This shows that case B3.3 is not realizable.

We note that items B3.1, B3.2, and B3.3 cover all the possibilities, so the consideration of the case B3 is complete.
This also completes the proof of the Case B and thus of Theorem~\ref{hugecases0} as well.
\end{proof}

The rest of this section is devoted to the consideration of the cases (i)--(v) of Theorem~\ref{hugecases0}.
In order to prove the main result of the section, we need to show that every of these cases is not realizable.

\subsection{Cases (ii) and (v)}

We start our consideration with the cases (ii) and (v) of Theorem~\ref{hugecases0}.

\begin{thr}\label{case VI}
Case (v) of Theorem~\ref{hugecases0} is not realizable.
\end{thr}

\begin{proof}
1. Let a matrix $W$ realize case (v), then $\trop(W)=3$, $\Kap(W)>3$,
and up to p.r.c. it holds that
$$\P(W)=\left(\begin{array}{c|c|c}
\i\ldots\i&\i\ldots\i&0\ldots0\\
\i\ldots\i&\i\ldots\i&0\ldots0\\
\i\ldots\i&\i\ldots\i&0\ldots0\\\hline
0\ldots0&\i\ldots\i&\i\ldots\i\\
0\ldots0&0\ldots0&\i\ldots\i\\
\underbrace{0\ldots0}_{s_1}&\underbrace{0\ldots0}_{s_2}&\underbrace{\i\ldots\i}_{s_3}\end{array}\right),\mbox{ where $s_1,s_2,s_3>0$}.$$
By Lemma~\ref{multiply}, w.l.o.g. we assume that the minimal element of every column of $W$ equals $0$.

2. Set $\ow=W[1,2,3,4,5]$. Then, by Lemma~\ref{monot}, $\trop(\ow)\leq3$. Now Theorem~\ref{cjrthr} shows that
$\Kap(\ow)\leq3$. By Definition~\ref{Kapr}, there exists a matrix $F'\in\Kf^{5\times n}$ such that
$\ow=\deg F'$ and $rank(F')\leq3$.

3. Thus every four rows of $F'$ are linearly dependent over $\Kf$. In particular, for any $k\in\{1,\ldots,n\}$ it holds that
\begin{equation}\label{VI1}\lambda_1 f'_{1k}+\lambda_2 f'_{2k}+\lambda_4 f'_{4k}+\lambda_5 f'_{5k}=0,
\mu_1 f'_{1k}+\mu_3 f'_{3k}+\mu_4 f'_{4k}+\mu_5 f'_{5k}=0.\end{equation}
We also set $\lambda_{3}=\lambda_{6}=\mu_{2}=\mu_{6}=0$.
Multiplying the linear combinations~(\ref{VI1}) by nonzero elements from $\Kf$, we assume w.l.o.g.
that $\min_{\theta=1}^6\{\deg\lambda_\theta\}=\min_{\theta=1}^6\{\deg\mu_\theta\}=0$. By Lemma~\ref{li->tr},
the tuples $(\deg\lambda_1 \deg \lambda_2 \deg \lambda_4 \deg \lambda_5)$ and $(\deg\mu_1 \deg\mu_3 \deg\mu_4 \deg\mu_5)$
realize the tropical dependence of rows of $W[1,2,4,5]$ and $W[1,3,4,5]$, respectively.

4. Lemma~\ref{r->b} shows that $\deg\lambda_1=\deg\lambda_2=\deg\mu_1=\deg\mu_3=0$, and the numbers
$\deg \lambda_4$, $\deg \mu_4$, $\deg \lambda_5$, $\deg \mu_5$ are all positive.

5. Since $\trop(W)\leq3$, there exists a tuple $(h_1 h_4 h_5 h_6)\in\T^4$ that realizes the tropical dependence of rows of $W[1,4,5,6]$.
Then by Lemma~\ref{r->b}, we assume w.l.o.g. that $h_1>0$, $h_4\geq0$, $h_5=h_6=0$. Therefore for any $\tau\in\{1,\ldots,s_1+s_2\}$
it holds that 
$h_5+\deg(f'_{5\tau})=0$, $h_4+\deg(f'_{4\tau})\geq0$.
The infiniteness of $\C$ allows us to find elements $\nu_4,\nu_5\in\Kf$, $\deg\nu_4=h_4$, $\deg\nu_5=h_5$, such that
\begin{equation}\label{VI2}deg(\nu_4 f'_{4\tau}+\nu_5 f'_{5\tau})=0 \mbox{ for any $\tau\in\{1,\ldots,s_1+s_2\}$}.\end{equation}
We also define $\nu_1$ and $\nu_6$ to be arbitrary elements of $\Kf$ of degrees $h_1$ and $h_6$, respectively, and set $\nu_2=\nu_3=0$.

6. Set $f_{i\tau}=f'_{i\tau}$ for $i\in\{1,2,3,4,5\}$, $\tau\in\{1,\ldots,s_1+s_2\}$. We also set
\begin{equation}\label{VI3}f_{6\tau}=\frac{-\nu_1 f_{1\tau}-\nu_4 f_{4\tau}-\nu_5 f_{5\tau}}{\nu_6}
\mbox{ for any $\tau\in\{1,\ldots,s_1+s_2\}$}.\end{equation}
Now item 5 shows that $\deg(\nu_1f_{1\tau})>0$, thus by~(\ref{VI2}),
it holds that $\deg f_{6\tau}=0$ for any $\tau\in\{1,\ldots,s_1+s_2\}$.

7. By item 5, the tuple $(\deg\nu_1 \deg\nu_4 \deg\nu_5 \deg\nu_6)$
realizes the tropical dependence of rows of $W[1,4,5,6]$. Thus by Lemma~\ref{sys1e},
for every $q\in\{s_1+s_2+1,\ldots,n\}$ there exist elements $f_{1q}, f_{4q}, f_{5q}, f_{6q}\in\Kf$ of the degrees
$w_{1q}, w_{4q}, w_{5q}, w_{6q}$, respectively, such that
\begin{equation}\label{VI4}\nu_1 f_{1q}+\nu_4 f_{4q}+\nu_5 f_{5q}+\nu_6 f_{6q}=0, \mbox{ where $q\in\{s_1+s_2+1,\ldots,n\}$}.\end{equation}

8. Finally, for any $q\in\{s_1+s_2+1,\ldots,n\}$ we set
\begin{equation}\label{VI5}f_{2q}=\frac{-\lambda_1 f_{1q}-\lambda_4 f_{4q}-\lambda_5 f_{5q}}{\lambda_2},
f_{3q}=\frac{-\mu_1 f_{1q}-\mu_4 f_{4q}-\mu_5 f_{5q}}{\mu_3}.\end{equation}
By item 4, $\deg\lambda_4>0$, $\deg\lambda_5>0$, so $\deg(\lambda_4f_{4q}+\lambda_5f_{5q})>0$, thus $\deg f_{2q}=0$
for any $q\in\{s_1+s_2+1,\ldots,n\}$. It can be shown in the same way that $\deg f_{3q}=0$.

9. Now the matrix $F\in\Kf^{6\times n}$ is well defined. Items 6, 7 and 8 show that $W=\deg F$.
By~(\ref{VI1}) and~(\ref{VI5}), both $2$nd and $3$rd rows of $F$ are linear combinations
of the $1$st, $4$th, and $5$th rows. By~(\ref{VI3}) and~(\ref{VI4}), the $6$th row of $F$ is also a linear combination
of the $1$st, $4$th, and $5$th rows. This shows that $rank(F)\leq 3$. By Definition~\ref{Kapr}, $\Kap(W)\leq3$.
The contradiction with item 1 shows that no such $W$ exists.
\end{proof}

Now we turn our attention to the case (ii) of Theorem~\ref{hugecases0}.

\begin{thr}\label{caseII}
Case (ii) of Theorem~\ref{hugecases0} is not realizable.
\end{thr}

\begin{proof}
1. Let a matrix $W$ realize case (ii), then $\trop(W)=3$, $\Kap(W)>3$, and
\begin{equation}\label{II1}\P(W)=\left(\begin{array}{c|c|c|c|c|c}
0\ldots0&\i\ldots\i&\i\ldots\i&\i\ldots\i&\i\ldots\i&\i\ldots\i\\
0\ldots0&\i\ldots\i&\i\ldots\i&\i\ldots\i&\i\ldots\i&\i\ldots\i\\\hline
\i\ldots\i&0\ldots0&0\ldots0&0\ldots0&\i\ldots\i&0\ldots0\\\hline
\i\ldots\i&0\ldots0&0\ldots0&\i\ldots\i&0\ldots0&0\ldots0\\\hline
\i\ldots\i&0\ldots0&\i\ldots\i&0\ldots0&0\ldots0&0\ldots0\\\hline
\underbrace{\i\ldots\i}_{s_0}&\underbrace{\i\ldots\i}_{s_1}&\underbrace{0\ldots0}_{s_2}&
\underbrace{0\ldots0}_{s_3}&\underbrace{0\ldots0}_{s_4}&\underbrace{0\ldots0}_{s_5}\end{array}\right),\end{equation}
we assume w.l.o.g. that $s_0>0$, $s_1>0$, $s_2+s_3+s_4+s_5>0$.

2. We assume w.l.o.g. that $s_1+\ldots+s_5$ is minimal over all matrices satisfying
the conditions of case (ii) of Theorem~\ref{hugecases0}.
By Lemma~\ref{multiply}, w.l.o.g. we assume that the minimal element of every column of $W$ equals $0$.

3. Then the minimal element of the matrix $W[1,2|s_0+1,\ldots,n]$ is $m>0$.
We add $-m$ to every element of the first two rows of $W$, $m$ to every element of the first $s_0$ columns of $W$.
The matrix $V$ obtained is such that
\begin{equation}\label{II2}\P(V)=\left(\begin{array}{c|c|c|c|c|c}
0\ldots0&P'_1&P'_2&P'_3&P'_4&P'_5\\
0\ldots0&&&&&\\\hline
\i\ldots\i&0\ldots0&0\ldots0&0\ldots0&\i\ldots\i&0\ldots0\\\hline
\i\ldots\i&0\ldots0&0\ldots0&\i\ldots\i&0\ldots0&0\ldots0\\\hline
\i\ldots\i&0\ldots0&\i\ldots\i&0\ldots0&0\ldots0&0\ldots0\\\hline
\underbrace{\i\ldots\i}_{s_0}&\underbrace{\i\ldots\i}_{s_1}&\underbrace{0\ldots0}_{s_2}&
\underbrace{0\ldots0}_{s_3}&\underbrace{0\ldots0}_{s_4}&\underbrace{0\ldots0}_{s_5}\end{array}\right),\end{equation}
and the matrix $\left(P'_1|\ldots|P'_5\right)$ contains at least one $0$. From item 2 it also follows that
the minimal element of every column of $V$ equals $0$.

4. Let us assume that the column $\left(\begin{smallmatrix}\i\\\i\end{smallmatrix}\right)$
appears in at least two of matrices $P'_1,\ldots,P'_5$.
The corresponding numbers of columns of $V$ are denoted by $\tau_1$ and $\tau_2$.
We assume w.l.o.g. that $\tau_1\in\{s_0+1,\ldots,s_0+s_1\}$.
Then by~(\ref{II2}), $[\P(V)]_{5\tau_1}=0$, $[\P(V)]_{6\tau_1}=\i$, $[\P(V)]_{6\tau_2}=0$.
If the first two rows of $\P(V)$ coincide,
then we add a small enough $\varepsilon>0$ to every element of the last four rows of $V$
and obtain a contradiction with the minimality of $s_1+\ldots+s_5$, assumed by item 2.

If otherwise the first two rows of $\P(V)$ are different, then $v_{1\tau_0}\neq v_{2\tau_0}$ for some $\tau_0$.
W.l.o.g. we assume $v_{1\tau_0}>0$, $v_{2\tau_0}=0$.
Then the matrix
$$\P(V)[5,6,2,1|\tau_1,\tau_2,\tau_0,1]=\left(
\begin{array}{cccc}
0&\pi_1&\pi_2&\i\\
\i&0&\pi_3&\i\\
\i&\i&0&0\\
\i&\i&\i&0\\
\end{array}
\right)
$$ is
$\B$-tropically singular,
so Lemma~\ref{Btrop} implies that $\trop(V)\geq4$.

The contradiction obtained shows that $\left(\begin{smallmatrix}\i\\\i\end{smallmatrix}\right)$
appears as a column of at most one of the matrices $P'_1,\ldots,P'_5$.
In particular, we assume w.l.o.g. that $\left(\begin{smallmatrix}\i\\\i\end{smallmatrix}\right)$ appears as a column of
none of $P'_2,P'_3,P'_4$.

5. By Lemma~\ref{monot}, $\trop(V[2,3,4,5,6])\leq3$.
Now Theorem~\ref{cjrthr} shows that $\Kap(V[2,3,4,5,6])\leq3$.
By Definition~\ref{Kapr}, there exists a matrix $F'$ such that
$V[2,3,4,5,6]=\deg F'$ and $rank(F')\leq3$.

6. By item 5, every four columns of $F'$ are linearly dependent over $\Kf$.
In particular, for every $j\in\{1,\ldots,n\}$ it holds that
\begin{equation}\label{II3}\lambda_{21} f'_{2j}+\lambda_{41} f'_{4j}+\lambda_{51} f'_{5j}+\lambda_{61} f'_{6j}=0,
\lambda_{22} f'_{2j}+\lambda_{32} f'_{3j}+\lambda_{42} f'_{4j}+\lambda_{62} f'_{6j}=0.\end{equation}
We set also $\lambda_{11}=\lambda_{31}=\lambda_{12}=\lambda_{52}=0$.

7. Multiplying the equations~(\ref{II3}) by elements from $\Kf^*$,
we assume w.l.o.g. that $\min_{k=1}^6\{\deg\lambda_{k1}\}=\min_{k=1}^6\{\deg\lambda_{k2}\}=0$.
By Lemma~\ref{li->tr}, the tuples $(\deg\lambda_{21} \deg \lambda_{41} \deg \lambda_{51} \deg \lambda_{61})$ and
$(\deg\lambda_{22} \deg\lambda_{32} \deg\lambda_{42} \deg\lambda_{62})$ realize
the tropical dependence of rows of $W[2,4,5,6]$ and $W[2,3,4,6]$, respectively.

8. By item 1, $s_1>0$, thus Lemma~\ref{r->b} implies that $\deg\lambda_{41}=\deg\lambda_{51}=\deg\lambda_{32}=\deg\lambda_{42}=0$, and the numbers
$\deg \lambda_{21}$, $\deg \lambda_{22}$ are positive. Analogously, if $s_3+s_4\neq0$, then $\deg\lambda_{62}=0$;
if $s_2+s_3\neq0$, then $\deg\lambda_{61}=0$.

9. Since $\trop(V)\leq3$, there is a tuple $(h_1 h_2 h_5 h_6)\in\T^4$, $\min\{h_1,h_2,h_5,h_6\}=0$, that realizes the tropical dependence
of rows of $W[1,2,5,6]$. Then by Lemma~\ref{r->b}, $h_1=h_2=0$. We set $\lambda_{13}$, $\lambda_{23}$, and $\lambda_{53}$
to be arbitrary elements from $\Kf$ of the degrees $h_1$, $h_2$, and $h_5$, respectively. Set also $\lambda_{33}=\lambda_{43}=0$.

10. Further, we denote $\rho_1=\lambda_{61}\lambda_{53}$, $\rho_2=\lambda_{51}\lambda_{63}$,
$\rho_3=(\lambda_{42}\lambda_{61}-\lambda_{62}\lambda_{41})\lambda_{53}$, $\rho_4=\lambda_{42}\lambda_{51}\lambda_{63}$.
From the infiniteness of $\C$ it follows that there exists $\lambda_{63}\in\Kf$ such that $\deg\lambda_{63}=h_6$,
$\deg(\rho_1-\rho_2)=\min\{\deg\rho_1,\deg\rho_2\}$, $\deg(\rho_3-\rho_4)=\min\{\deg\rho_3,\deg\rho_4\}$.
So we have defined the elements $\{\lambda_{ki}\}$ for every $k\in\{1,\ldots,6\}$, $i\in\{1,2,3\}$.
The matrix $(\lambda_{ki})$ is further denoted by $\Lambda\in\Kf^{6\times3}$.

11. If $s_2\neq0$, then from item 8 and~(\ref{II3}) it follows that
$\deg(\lambda_{41}f'_{4,s_1+1}+\lambda_{61}f'_{6,s_1+1})>0$, $\deg(\lambda_{42}f'_{4,s_1+1}+\lambda_{62}f'_{6,s_1+1})=0$.
Thus if $s_2\neq0$, then $\deg(\lambda_{41}\lambda_{62}-\lambda_{42}\lambda_{61})=0$. In the same way, we can prove that
if $s_4\neq0$, then also $\deg(\lambda_{41}\lambda_{62}-\lambda_{42}\lambda_{61})=0$.

12. We start the construction of a matrix
$F$$\in$$\Kf^{6\times n}$ such that $\deg$$F$$=$$V$ and
\begin{equation}\label{II51}\lambda_{1i}f_{1j}+\lambda_{2i}f_{2j}+\lambda_{3i}f_{3j}+\lambda_{4i}f_{4j}+\lambda_{5i}f_{5j}+\lambda_{6i}f_{6j}=0\end{equation}
for any $i\in\{1,2,3\}$, $j\in\{1,\ldots,n\}$. We consider the five cases.

\textbf{Case A.} Let us assume that $j\in\{1,\ldots,n\}$ is such that the minimum in
\begin{equation}\label{II52}\min_{\theta=1}^6\{\deg\lambda_{\theta3}+v_{\theta j}\}\end{equation}
is attained if and only if $\theta\in\{1,2\}$. From item 9 it then follows that $v_{1j}=v_{2j}$.
Set $f_{kj}=f'_{kj}$ for $k\neq1$, set also
\begin{equation}\label{II53}f_{1j}=\frac{-\lambda_{23}f_{2j}-\lambda_{53}f_{5j}-\lambda_{63}f_{6j}}{\lambda_{13}}.\end{equation}
The only term of degree $v_{2j}$ in the numerator is $\lambda_{23}f_{2j}$, all the other terms have greater degrees,
so $\deg f_{1j}=v_{2j}=v_{1j}$. Hence from item 5 it follows that $\deg f_{kj}=v_{kj}$ for every $k\in\{1,2,3,4,5,6\}$.
The equations~(\ref{II51}) for $i=1$ and $i=2$ now follow from item 6, for $i=3$ from~(\ref{II53}).

\textbf{Case B.} Now we assume $j\in\{s_0+1,\ldots,s_0+s_1,s_0+s_1+s_2+s_3+s_4+1,\ldots,n\}$ and consider the two possible cases.

Case B1. Let $\theta=5$ and $\theta=6$ provide the minimum for~(\ref{II52}).
Then $v_{6j}=h_5-h_6$ and $\min_{\theta=1}^6\{\deg\lambda_{\theta 3}+v_{\theta j}\}=h_5$.
Now items 8 and 9 imply that $\deg\det\Lambda[3,4,5]=h_5$, $\deg\det\Lambda[3,4,6]=h_6$,
and for any $i'\in\{1,2\}$ it holds that $\min_{\theta=1}^6\{\deg\lambda_{\theta i'}+v_{\theta j}\}=0$.
From item 10 it also follows that $\deg\det\Lambda[3,5,6]=\deg\det\Lambda[4,5,6]=h_6$.
Thus we see that the matrix $\Lambda\in\Kf^{6\times3}$,
the tuple $(v_{1j},\ldots,v_{6j})$, and the indexes $\{3,4,5,6\}$ satisfy the conditions of Lemma~\ref{sys}.
This implies that for some $f_{1j},\ldots,f_{6j}\in\Kf$ it holds that $\deg f_{kj}=v_{kj}$ for any $k\in\{1,2,3,4,5,6\}$,
and the equations~(\ref{II51}) hold for $i\in\{1,2,3\}$.

Case B2. Let either $\theta=1$ or $\theta=2$ provide the minimum for~(\ref{II52}).
We set $f_{6j}=t^{v_{6j}}$. By Lemma~\ref{forsystem2}, there exists $\xi\in\C^*$ such that
\begin{equation}\label{CaseB2ii}
\deg\left(\lambda_{53}\xi+\lambda_{63}f_{6j}\right)=\min\{\deg\lambda_{53},\deg\lambda_{63}+v_{6j}\},
\end{equation}
\begin{equation}\label{CaseB2ii2}
\deg\left(\lambda_{51}\xi+\lambda_{61}f_{6j}\right)=0,
\end{equation}
\begin{equation}\label{CaseB2ii3}
\deg\left(\frac{\lambda_{42}\lambda_{51}}{\lambda_{41}}\xi+\left(\frac{\lambda_{42}\lambda_{61}}{\lambda_{41}}-\lambda_{62}\right)f_{6j}\right)=0.
\end{equation}
We set $f_{5j}=\xi$. We see that $\deg f_{5j}=v_{5j}$, $\deg f_{6j}=v_{6j}$. Since the minimum in~(\ref{II52})
is provided by either $\theta=1$ or $\theta=2$, from~(\ref{CaseB2ii}) it follows that there exist
elements $f_{1j},f_{2j}\in\Kf$ such that $\deg f_{1j}=v_{1j}$, $\deg f_{2j}=v_{2j}$,
$\lambda_{13}f_{1j}+\lambda_{23}f_{2j}+\lambda_{53}f_{5j}+\lambda_{63}f_{6j}=0$.
In this case, the condition~(\ref{II51}) holds for $i=3$ for any $f_{3j}$ and $f_{4j}$
because item 9 implies that $\lambda_{33}=\lambda_{43}=0$.

Further, we set
\begin{equation}\label{II54}f_{4j}=\frac{-\lambda_{21}f_{2j}-\lambda_{51}f_{5j}-\lambda_{61}f_{6j}}{\lambda_{41}},
f_{3j}=\frac{-\lambda_{22}f_{2j}-\lambda_{42}f_{4j}-\lambda_{62}f_{6j}}{\lambda_{32}}.\end{equation}
By item 6, $\lambda_{11}=\lambda_{31}=\lambda_{12}=\lambda_{52}=0$,
thus the conditions~(\ref{II51}) with $i=1$ and $i=2$ follow from~(\ref{II54}).

Finally, item 8 implies that $\deg(\lambda_{21}f_{2j})>0$, so from~(\ref{CaseB2ii2}) it follows that $f_{4j}=0=v_{4j}$.
In the definition of $f_{3j}$, we substitute the values of $f_{4j}$ and $f_{5j}$ by their expressions and obtain
$$f_{3j}=\frac{\lambda_{42}\lambda_{51}}{\lambda_{32}\lambda_{41}}\xi+\left(\frac{\lambda_{61}\lambda_{42}}{\lambda_{32}\lambda_{41}}-
\frac{\lambda_{62}}{\lambda_{32}}\right)f_{6j}+\left(\frac{\lambda_{21}\lambda_{42}}{\lambda_{32}\lambda_{41}}-\frac{\lambda_{22}}{\lambda_{32}}\right)f_{2j}.$$
From item 8 and~(\ref{CaseB2ii3}) it follows that $f_{3j}=0=v_{3j}$. This completes the consideration of Case B2.

Item 9 shows that $(\deg\lambda_{13},\deg\lambda_{23},\deg\lambda_{53},\deg\lambda_{63})$ realizes the tropical dependence of rows of $V[1,2,5,6]$,
so Cases B1 and B2 cover all the possibilities. The consideration of Case B is complete.

\textbf{Case C.} Now we assume that $j\in\{s_0+s_1+1,\ldots,s_0+s_1+s_2\}$, and that the assumption of Case A fails to hold for $j$.

C1. Item 4 shows that $v_{g_1j}=0$ for some $g_1\in\{1,2\}$, item 9 now implies that $\deg\lambda_{g_13}=0$.
Thus the minimum in~(\ref{II52}) is attained for $\theta=g_1$.
By item 9, the tuple $(\deg\lambda_{13},\deg\lambda_{23},\deg\lambda_{53},\deg\lambda_{63})$
realizes the tropical dependence of rows of $V[1,2,5,6]$, thus the minimum in~(\ref{II52}) is also attained
for some $g_2\neq g_1$. Since the assumption of Case A fails to hold for $j$, we assume w.l.o.g. that
$g_2\notin\{1,2\}$. From the equation~(\ref{II2}) it follows that $g_2\neq5$, from item 9 that $g_2\notin\{3,4\}$, so $g_2=6$.
Then $\deg\lambda_{63}+v_{6j}=\deg\lambda_{g_13}+v_{g_1j}=0$, i.e. $\deg\lambda_{63}=0$.

C2. Now from items 8 and 9 it follows that $\deg\det\Lambda[g_1,3,4]=\deg\det\Lambda[g_1,3,6]=0$,
and for any $i\in\{1,2,3\}$ it holds that $\min_{\theta=1}^6\{\deg\lambda_{\theta i}+v_{\theta j}\}=0$.
The equality $\deg\det\Lambda[g_1,4,6]=0$ follows from item 11, $\deg\det\Lambda[3,4,6]=0$ from items 8 and C1.
Thus we see that the matrix $\Lambda\in\Kf^{6\times3}$,
the tuple $(v_{1j},\ldots,v_{6j})$, and the indexes $\{g_1,3,4,6\}$ satisfy the conditions of Lemma~\ref{sys}.
This implies that for some $f_{1j},\ldots,f_{6j}\in\Kf$ it holds that $\deg f_{kj}=v_{kj}$ for any $k\in\{1,2,3,4,5,6\}$,
and the equations~(\ref{II51}) hold for $i\in\{1,2,3\}$.

\textbf{Case D.} Let us now assume that $j\in\{s_0+s_1+s_2+1,\ldots,s_0+s_1+s_2+s_3\}$, and that the assumption of Case A fails to hold for $j$.
The argument similar to one of item C1 shows that the minimum in~(\ref{II52}) is then attained
for some $\theta_1\in\{1,2\}$ and $\theta_2\in\{5,6\}$, we also obtain that $\deg\lambda_{\theta_23}=0$.
Now items 8 and 9 imply that $\deg\det\Lambda[\theta_1,3,5]=\deg\det\Lambda[\theta_1,3,6]=\deg\det\Lambda[\theta_1,5,6]=0$,
and for any $i\in\{1,2,3\}$ it holds that $\min_{\theta=1}^6\{\deg\lambda_{\theta i}+v_{\theta j}\}=0$.
The equality $\deg\det\Lambda[3,5,6]=0$ follows from item 10.
Thus we see that the matrix $\Lambda\in\Kf^{6\times3}$,
the tuple $(v_{1j},\ldots,v_{6j})$, and the indexes $\{\theta_1,3,5,6\}$ satisfy the conditions of Lemma~\ref{sys}.
This implies that for some $f_{1j},\ldots,f_{6j}\in\Kf$ it holds that $\deg f_{kj}=v_{kj}$ for any $k\in\{1,2,3,4,5,6\}$,
and the equations~(\ref{II51}) hold for $i\in\{1,2,3\}$.

\textbf{Case E.} Now assume that $j\in\{s_0+s_1+s_2+s_3+1,\ldots,s_0+s_1+s_2+s_3+s_4\}$, and that the assumption of Case A fails to hold for $j$.
The argument similar to one of item C1 shows that the minimum in~(\ref{II52}) is then attained
for some $\theta'_1\in\{1,2\}$ and $\theta'_2\in\{5,6\}$.
Now item 8 implies that $\deg\det\Lambda[\theta'_1,4,5]=\deg\det\Lambda[\theta'_1,5,6]=0$,
and for any $i\in\{1,2,3\}$ it holds that $\min_{\theta=1}^6\{\deg\lambda_{\theta i}+v_{\theta j}\}=0$.
From items 8 and 11 it also follows that $\deg\det\Lambda[\theta'_1,4,6]=0$, from item 10 that $\deg\det\Lambda[4,5,6]=0$.
Thus we see that the matrix $\Lambda\in\Kf^{6\times3}$,
the tuple $(v_{1j},\ldots,v_{6j})$, and the indexes $\{\theta'_1,4,5,6\}$ satisfy the conditions of Lemma~\ref{sys}.
This implies that for some $f_{1j},\ldots,f_{6j}\in\Kf$ it holds that $\deg f_{kj}=v_{kj}$ for any $k\in\{1,2,3,4,5,6\}$,
and the equations~(\ref{II51}) hold for $i\in\{1,2,3\}$.

Now we note that Cases A--E cover all the possibilities for the number of column $j$.
Indeed, item 9 implies that the numbers $j'\in\{1,\ldots,s_0\}$ satisfy the assumption of Case A.
All the other possibilities have been considered in Cases B--E.
Thus we see that there exists a matrix $F$ such that $V=\deg F$, and the conditions~(\ref{II51}) hold
for every $i\in\{1,2,3\}$ and $j\in\{1,\ldots,n\}$. By the construction of $\Lambda$, this implies that
every row of $F$ is a linear combination of its $2$nd, $4$th, and $6$th rows. This shows that $rank(F)\leq 3$.
By Definition~\ref{Kapr}, $\Kap(W)\leq3$. The contradiction with item 1 shows that no such $W$ exists.
\end{proof}

\subsection{Case (iv)}

This subsection deals with case (iv) of Theorem~\ref{hugecases0}.
We start with the following lemmas.

\begin{lem}\label{startV310711}
Let the tropical rank of a matrix $M\in\R^{4\times n}$ be at most $3$.
Let $\left(\begin{smallmatrix}0\\\i\\\i\\\i\end{smallmatrix}\right)$ appear as a column of $\P(M)$,
the rows of $\P(M)[2,3,4]$ be pairwise different, every row of $\P(M)$ contain at least one $0$.
Then there exists a positive number $x$ such that the tuple $(x,0,0,0)$ realizes
the tropical dependence of rows of $M$.
\end{lem}

\begin{proof}
By Theorem~\ref{tropthr}, there exists a tuple $(x,a,b,c)$ realizing the tropical dependence of rows of $M$.
We assume w.l.o.g. that $\min\{x,a,b,c\}=0$. From Lemma~\ref{r->b} it then follows that $x>0$ and $a=b=c=0$.
\end{proof}

\begin{lem}\label{startV}
Let a matrix $W$ realize case (iv) of Theorem~\ref{hugecases0}.
Then there exist positive numbers $x$ and $y$ such that the tuples  $(0,0,0,x)$ and $(y,0,0,0)$ realize
the tropical dependence of rows of $W[1,2,3,4]$ and $W[3,4,5,6]$, respectively.
\end{lem}

\begin{proof}
1. Under the assumptions of (iv), the matrix $W$ is such that $\trop(W)=3$, $\Kap(W)>3$, and
$$\P(W)=\left(\begin{array}{c|c|c}
0\ldots0&\i\ldots\i&\\
0\ldots0&\i\ldots\i&P'\\
0\ldots0&\i\ldots\i&\\\hline
\i\ldots\i&0\ldots0&\\
\i\ldots\i&0\ldots0&P''\\
\underbrace{\i\ldots\i}_{p}&\underbrace{0\ldots0}_{q}&\underbrace{}_{r}\end{array}\right),$$
where $p$ and $q$ are nonzero, $r$ may equal zero, every column of $P'$ and every column of $P''$ contain
at least two zeros. By Lemma~\ref{multiply}, w.l.o.g. we assume that the minimal element of
every column of $W$ equals $0$.

2. We will only prove that $(y,0,0,0)$ realizes the tropical dependence of $W[3,4,5,6]$ for some $y>0$.
The case of $W[1,2,3,4]$ can be considered in the same way. The two cases are possible.

\textbf{Case A.} Let the element $\i$ appear as an entry of $P''$. In this case, if some rows of $P''$ coincide,
then we add a small enough $-\varepsilon<0$ to every element of the last three rows of $W$. By Definition~\ref{patt},
the matrix obtained satisfies up to p.r.c. the conditions of case (v) of Theorem~\ref{hugecases0}.

The contradiction with Theorem~\ref{case VI} shows that the rows of $P''$ are pairwise different.
Then Lemma~\ref{startV310711} completes the consideration of Case A.

\textbf{Case B.} Let the matrix $P''$ consist of zero elements. This case is treated by
\textit{reductio ad absurdum}. We assume that for every $y>0$ the tuple $(y,0,0,0)$ does
not realize the tropical dependence of rows of $W[3,4,5,6]$.

B1. Then, by Definition~\ref{tropdep}, there exist $j_1,j_2\in\{1,\ldots,p\}$ such that the minimum
over the set $\{w_{4j_1},w_{5j_1},w_{6j_1},w_{4j_2},w_{5j_2},w_{6j_2}\}$ is attained exactly once.
We assume w.l.o.g. that $a=w_{4j_1}<\min\{w_{5j_1},w_{6j_1},w_{4j_2},w_{5j_2},w_{6j_2}\}$.

B2. By $J'$ we denote the set of all $j$ such that $w_{4j}, w_{5j}, w_{6j}$ are not all equal.
We set $$a'=\min_{i\in\{4,5,6\}, j\in J'}\{w_{ij}\}.$$ By item B1, $j_1\in J'$ and $a'\leq a$.

B3. We now add $-a'$ to every element of the last three rows of $W$.
By Definition~\ref{patt}, the pattern of the matrix $V$ obtained is
\begin{equation}\label{V01}\P(V)=\left(\begin{array}{c|c}
0\ldots0&\\
0\ldots0&\widehat{P}\\
0\ldots0&\\\hline
&0\ldots0\\
\widetilde{P}&0\ldots0\\
&0\ldots0\end{array}\right),\end{equation}
where $\widetilde{P}$ contains at least one $0$, and $\left(\begin{smallmatrix}0\\0\\0\end{smallmatrix}\right)$
does not appear as a column of $\widetilde{P}$. Note that by item B1, the $j_2$-th column of $\widetilde{P}$
is $\left(\begin{smallmatrix}\i\\\i\\\i\end{smallmatrix}\right)$. We consider the two possible cases.

B4. Assume that the last two rows of $\widetilde{P}$ are different. Then w.l.o.g. we assume that
$[\P(V)]_{5j_0}=0$, $[\P(V)]_{6j_0}=\i$ for some $j_0$. This implies that $v_{6j_0}>v_{5j_0}=0$.
By item B3, we have
$$V[3,4,5,6|j_2,j_1,j_0,p'+1]=
\left(
\begin{array}{cccc}
0&0&0&z_1\\
w_{4j_2}-a'&a-a'&z_2&0\\
w_{5j_2}-a'&w_{5j_1}-a'&0&0\\
w_{6j_2}-a'&w_{6j_1}-a'&v_{6j_0}&0
\end{array}
\right),
$$
where $z_1,z_2$ are nonnegative numbers. Item B1 implies that $a-a'<\min\{w_{5j_1}-a',w_{6j_1}-a',w_{4j_2}-a',w_{5j_2}-a',w_{6j_2}-a'\}$.
Thus by Definition~\ref{tropdeg}, $V[3,4,5,6|j_2,j_1,j_0,p'+1]$ is tropically nonsingular. Thus $\trop(W)\geq4$,
so we have a contradiction with item 1.

B5. So item B4 shows that the last two rows of $\widetilde{P}$ coincide.

B5.1. In this case, if the first row of $\widetilde{P}$ contains at least one $0$, then
we add a small enough $\varepsilon>0$ to every element of the first four rows of $V$.
By Definition~\ref{patt}, the matrix obtained up to p.r.c. satisfies the conditions of case (ii)
of Theorem~\ref{hugecases0}. This contradicts Theorem~\ref{caseII}.

B5.2. If otherwise the elements of the first row of $\widetilde{P}$ are all equal to $\i$,
then we add a small enough $-\varepsilon<0$ to every element of the last three rows of $V$.
By Definition~\ref{patt}, the matrix obtained up to p.r.c. satisfies the conditions of case (v)
of Theorem~\ref{hugecases0}. This contradicts Theorem~\ref{case VI}.

The contradiction obtained completes the consideration of Case B, showing that the tuple $(y,0,0,0)$
realizes the tropical dependence of rows of $W[3,4,5,6]$ for some $y>0$.
\end{proof}

We are now ready to prove the main result of this subsection.

\begin{thr}\label{cases IV-V}
Case (iv) of Theorem~\ref{hugecases0} is not realizable.
\end{thr}

\begin{proof}
1. Let a matrix $W$ realize case (iv), then $\trop(W)=3$, $\Kap(W)>3$.
Also, then by Lemma~\ref{multiply}, we can assume w.l.o.g. that $W$ is nonnegative and
\begin{equation}\label{V1}W=\left(\begin{array}{c|c|c}
0\ldots0&&\\
0\ldots0&\widetilde{P'}&P'\\
0\ldots0&&\\\hline
&0\ldots0&\\
\widetilde{P''}&0\ldots0&P''\\
\underbrace{}_{p}&\underbrace{0\ldots0}_{q}&\underbrace{}_{r}\end{array}\right),\end{equation}
where the matrices $\widetilde{P'}$ and $\widetilde{P''}$ consist of positive elements,
every column of $P'$ and $P''$ contains at least two zeros.

2. By Lemma~\ref{startV}, there exist positive numbers $a$, $b$ such that the tuples $(0,0,0,a)$ and $(b,0,0,0)$
realize the tropical dependence of rows of $W[1,2,3,4]$ and $W[3,4,5,6]$, respectively.

3. By $\Lambda\in\Kf^{3\times3}$ we denote a matrix whose entries have the degree $0$,
cofactors the degree $a$, determinant the degree $2a+b$.
To be definite, set
$$\Lambda=\left(
\begin{array}{ccc}
1+t^a&1&1\\
1&1+t^a&1\\
1&1&\frac{2}{2+t^a}+t^{a+b}
\end{array}
\right).
$$

4. We will show that for every $j\in\{1,\ldots,n\}$ there exist elements $f_{1j},\ldots,f_{6j}$
of the degrees $w_{1j},\ldots,w_{6j}$, respectively, such that
\begin{equation}\label{sysIVV}
\Lambda\left(
\begin{array}{c}
f_{1j}\\
f_{2j}\\
f_{3j}\\
\end{array}
\right)=
\left(
\begin{array}{c}
t^af_{4j}\\
t^af_{5j}\\
t^af_{6j}\\
\end{array}
\right).
\end{equation}

We consider the possible cases.

\textbf{Case A.} Assume that $j\in\{1,\ldots,p\}$, and the minimum over the set $\{w_{4j},w_{5j},w_{6j}\}$ is attained exactly once.
From item 2 it now follows that $\min\{w_{4j},w_{5j},w_{6j}\}=b.$ Then for every $i'\in\{1,2,3\}$ the element
$$f_{i'j}=\frac{t^{a+w_{4j}}\Lambda_{1i'}+t^{a+w_{5j}}\Lambda_{2i'}+t^{a+w_{6j}}\Lambda_{3i'}}{\det\Lambda}$$
has a zero degree. We also set $f_{4j}=t^{w_{4j}}$, $f_{5j}=t^{w_{5j}}$, $f_{6j}=t^{w_{6j}}$.
The Cramer's rule for solving linear systems shows that $f_{1j},\ldots,f_{6j}$ satisfy the condition~(\ref{sysIVV}).

\textbf{Case B.} Assume that $j\in\{1,\ldots,p\}$, and the minimum over the set $\{w_{4j},w_{5j},w_{6j}\}$ is attained at least twice.
From item 2 it now follows that $\min\{w_{4j},w_{5j},w_{6j}\}=c\leq b.$ Equation~(\ref{V1}) implies that $c>0$.
We assume w.l.o.g. that $w_{4j}=w_{5j}=c$, $w_{6j}=d\geq c$.
Then the degree of \begin{equation}\label{f5j}f_{5j}=-\frac{\gamma_1 t^c\Lambda_{11}}{\Lambda_{21}}-\frac{t^d\Lambda_{31}}{\Lambda_{21}}+t^{b}\end{equation}
equals $c$ for some $\gamma_1\in\C^*$. We also set $f_{4j}=\gamma_1t^c$, $f_{6j}=t^d$.

Now for $i'\in\{1,2,3\}$ we set $$f_{i'j}=\frac{t^a f_{4j}\left(\Lambda_{1i'}-\frac{\Lambda_{2i'}\Lambda_{11}}{\Lambda_{21}}\right)+
t^a f_{6j}\left(\Lambda_{3i'}-\frac{\Lambda_{2i'}\Lambda_{31}}{\Lambda_{21}}\right)+t^{a+b}\Lambda_{2i'}}{\det\Lambda}.$$
Item 3 implies that $f_{1j}=f_{2j}=f_{3j}=0$. From~(\ref{f5j}) it follows that
$$f_{i'j}=\frac{t^af_{4j}\Lambda_{1i'}+t^af_{5j}\Lambda_{2i'}+t^af_{6j}\Lambda_{3i'}}{\det\Lambda},$$
so the Cramer's rule implies the condition~(\ref{sysIVV}).

\textbf{Case C.} Now assume that $j\in\{p+1,\ldots,p+q\}$, and the minimum over the set $\{w_{1j},w_{2j},w_{3j}\}$ is attained exactly once.
From item 2 it now follows that $\min\{w_{1j},w_{2j},w_{3j}\}=a.$
Then for every $i'\in\{1,2,3\}$ the element
$$f_{i'+3,j}=t^{-a}\left(\lambda_{i'1}t^{w_{1j}}+\lambda_{i'2}t^{w_{2j}}+\lambda_{i'3}t^{w_{3j}}\right)$$
has a zero degree. We set $f_{1j}=t^{w_{1j}}$, $f_{2j}=t^{w_{2j}}$, $f_{3j}=t^{w_{3j}}$.
These settings satisfy the condition~(\ref{sysIVV}).

\textbf{Case D.} Assume that $j\in\{p+1,\ldots,p+q\}$, and the minimum over the set $\{w_{1j},w_{2j},w_{3j}\}$ is attained at least twice.
From item 2 it now follows that $\min\{w_{1j},w_{2j},w_{3j}\}=h\leq a.$ Equation~(\ref{V1}) implies that $h>0$.
We assume w.l.o.g. that $w_{1j}=w_{2j}=h$, $w_{3j}=g\geq h$. Then the degree of
\begin{equation}\label{f2j}f_{2j}=-\frac{\gamma_2t^h\lambda_{11}}{\lambda_{12}}-\frac{t^g\lambda_{13}}{\lambda_{12}}+t^{a}\end{equation}
equals $h$ for some $\gamma_2\in\C^*$. We also set $f_{1j}=\gamma_2t^h$, $f_{3j}=t^g$.

Now for $i'\in\{1,2,3\}$ we set
$f_{i'+3,j}=t^{-a}\left(\lambda_{i'1}f_{1j}+\lambda_{i'2}f_{2j}+\lambda_{i'3}f_{3j}\right).$
These settings satisfy the condition~(\ref{sysIVV}).
From~(\ref{f2j}) it now follows that
$$f_{i'+3,j}=\frac{\gamma_2t^h\left(\lambda_{i'1}-\frac{\lambda_{11}\lambda_{i'2}}{\lambda_{12}}\right)+
t^g\left(\lambda_{i'3}-\frac{\lambda_{13}\lambda_{i'2}}{\lambda_{12}}\right)+\lambda_{i'2}t^{a}}{t^a},$$
so item 3 implies that $f_{4j}=f_{5j}=f_{6j}=0$.

\textbf{Case E.} Finally, let $j\in\{p+q+1,\ldots,n\}$. By item 1, it can be assumed w.l.o.g. that
$w_{1j}=w_{2j}=0=w_{4j}=w_{5j}=0$, $w_{3j}=\alpha\geq0$, $w_{6j}=\beta\geq0$.
Then all but one complex numbers $\zeta$ are such that the element
$$f_{4j}=-\frac{\zeta\Lambda_{23}}{\Lambda_{13}}-\frac{t^\beta\Lambda_{33}}{\Lambda_{13}}+t^{b+\alpha}$$
has a zero degree. We also assume $\zeta\neq0$ and set $f_{5j}=\zeta$, $f_{6j}=t^\beta$.
We also set $$f_{i'j}=\frac{t^af_{4j}\Lambda_{1i'}+t^af_{5j}\Lambda_{2i'}+t^af_{6j}\Lambda_{3i'}}{\det\Lambda}$$ for $i'\in\{1,2,3\}$,
then the Cramer's rule implies the condition~(\ref{sysIVV}).

Our settings imply that for $i'\in\{1,2,3\}$ it holds that
$$f_{i'j}=\frac{t^{a+\beta}\left(\Lambda_{3i'}-\frac{\Lambda_{33}\Lambda_{1i'}}{\Lambda_{13}}\right)+
\zeta t^{a}\left(\Lambda_{2i'}-\frac{\Lambda_{23}\Lambda_{1i'}}{\Lambda_{13}}\right)+t^{a+b+\alpha}\Lambda_{1i'}}{\det\Lambda},$$
in particular, $f_{3j}=\frac{t^{a+b+\alpha}\Lambda_{13}}{\det\Lambda}$. From item 3 it follows that $\deg f_{3j}=\alpha$, and
all but one or two complexes $\zeta$ are such that $\deg f_{1j}=\deg f_{2j}=0$. The infiniteness of $\C$ implies
that there exists $\zeta\in\C$ such that $\deg f_{ij}=w_{ij}$ for every $i\in\{1,2,3,4,5,6\}$.

Cases A--E cover all the possibilities, so there exists a lift $F$ of $W$ such that
the condition~(\ref{sysIVV}) is satisfied. Thus $rank F\leq3$, so by Definition~\ref{Kapr},
$\Kap(W)\leq3$. The contradiction with item 1 shows that no such $W$ exists.
\end{proof}

\subsection{Case (iii)}

This subsection is devoted to case (iii) of Theorem~\ref{hugecases0}.
We need the following lemma.

\begin{lem}\label{pre-iii}
Let a matrix $W$ realize case (iii) of Theorem~\ref{hugecases0},
$\{r_1,r_2,r_3\}$ be a support of some column of $\P(W)$, $\{r_4,r_5,r_6\}=\{1,\ldots,6\}\setminus\{r_1,r_2,r_3\}$.
Then the rows of $\P(W)[r_4,r_5,r_6]$ are pairwise distinct.
\end{lem}

\begin{proof}
Assume the converse. Then w.l.o.g. we can assume that $\{r_1,r_2,r_3\}=\{1,2,3\}$,
and that the fourth and fifth rows of $\P(W)$ coincide. Then
$$\P(W)=\left(\begin{array}{c|c|c|c}
0\ldots0&&&\\
0\ldots0&P_2&P_3&P_4\\
0\ldots0&&&\\\hline
\i\ldots\i&0\ldots0&0\ldots0&\i\ldots\i\\
\i\ldots\i&0\ldots0&0\ldots0&\i\ldots\i\\\hline
\underbrace{\i\ldots\i}_{h_1}&\underbrace{0\ldots0}_{h_2}&\underbrace{\i\ldots\i}_{h_3}&\underbrace{0\ldots0}_{h_4}\end{array}\right),$$
where $h_2$, $h_3$, and $h_4$ may equal $0$. The assumptions of case (iii) imply that indeed $h_4=0$. Moreover, Theorem~\ref{hugecases0} requires
every row of $\P(W)$ to contain at least one $0$, so we have $h_2\neq0$.

Now we add a small enough $\varepsilon>0$ to every element of the first three rows of $W$. By Definition~\ref{patt},
the matrix obtained satisfies up to p.r.c. the conditions of either case (v) (if $h_3\neq0$) or case (iv) (if $h_3=0$)
of Theorem~\ref{hugecases0}. The contradiction with Theorems~\ref{case VI} and~\ref{cases IV-V} completes the proof.
\end{proof}

\begin{thr}\label{caseIIId}
Case (iii) of Theorem~\ref{hugecases0} is not realizable.
\end{thr}

\begin{proof}
1. Let a matrix $W$ realize case (iii). Then, in particular, $\trop(W)=3$, $\Kap(W)>3$, and
\begin{equation}\label{eqIIId1}\P(W)=\left(\begin{array}{c|c|c|c|c}
0\ldots0&0\ldots0&\i\ldots\i&\i\ldots\i&\\
0\ldots0&\i\ldots\i&0\ldots0&\i\ldots\i&\\
0\ldots0&\i\ldots\i&\i\ldots\i&0\ldots0&\mbox{$ $ $P'$ $ $}\\
\i\ldots\i&0\ldots0&0\ldots0&\i\ldots\i&\\
\i\ldots\i&0\ldots0&\i\ldots\i&0\ldots0&\\
\underbrace{\i\ldots\i}_{q_1}&\underbrace{\i\ldots\i}_{q_2}&\underbrace{0\ldots0}_{q_3}&\underbrace{0\ldots0}_{q_4}&\end{array}\right),\end{equation}
where any of $q_1,\ldots,q_4$ may equal $0$.

2. By Lemma~\ref{multiply}, w.l.o.g. we assume that the minimal element of every column of $W$ equals $0$.

3. If $q_1\neq0$, then Lemmas~\ref{startV310711} and~\ref{pre-iii} imply that there exists $a>0$ such that the tuple
$(a,0,0,0)$ realizes the tropical dependence of rows of $W[3,4,5,6]$. If $q_1=0$, then we set $a=0$.

Analogously, if $q_2\neq0$, then there exists $b>0$ such that
$(b,0,0,0)$ realizes the tropical dependence of rows of $W[1,2,3,6]$. If $q_2=0$, then we set $b=0$.

If $q_3\neq0$, then there exists $c>0$ such that the tuple
$(0,c,0,0)$ realizes the tropical dependence of rows of $W[1,2,3,5]$. If $q_3=0$, then we set $c=0$.

If $q_4\neq0$, then there exists $d>0$ such that the tuple
$(0,0,d,0)$ realizes the tropical dependence of rows of $W[1,2,3,4]$. If $q_4=0$, then we set $d=0$.

4. We set
$$\Lambda=\left(
\begin{array}{ccc}
2+t^a&1&-t^b\\
1&2t^c&1\\
5t^d&1&\frac{2+t^b-10t^{b+c+d}-5t^d}{-1+4t^c+2t^{a+c}}\\
-1&0&0\\
0&-1&0\\
0&0&-1
\end{array}
\right).
$$
We also denote $L=\Lambda[1,2,3]$.

5. Let $g_1, g_2, g_3\in\{1,2,3,4,5,6\}$ be pairwise distinct numbers. One can check that $\deg\det\Lambda[g_1,g_2,g_3]=0$
if $\{g_1,g_2,g_3\}$ is a support of no column of $\P(W)$. Also, the computation shows that $\deg\det L=a$,
and every cofactor of $L$ has a zero degree.

6. We will show that for every $j\in\{1,\ldots,n\}$ there exist elements $f_{1j},\ldots,f_{6j}$
of the degrees $w_{1j},\ldots,w_{6j}$, respectively, such that
\begin{equation}\label{eqIIId108}\lambda_{1i}f_{1j}+\ldots+\lambda_{6i}f_{6j}=0 \mbox{ for every }i\in\{1,2,3\}.\end{equation}
We consider the four cases.

\textbf{Case A.} First, assume that $j\in\{1,\ldots,q_1\}$, and the minimum over the set $\{w_{4j},w_{5j},w_{6j}\}$ is attained exactly once.
From item 3 it now follows that $\min\{w_{4j},w_{5j},w_{6j}\}=a.$ Then by item 5, for every $i'\in\{1,2,3\}$ the element
$$f_{i'j}=\frac{t^{w_{4j}}L_{i'1}+t^{w_{5j}}L_{i'2}+t^{w_{6j}}L_{i'3}}{\det L}$$
has a zero degree. We also set $f_{4j}=t^{w_{4j}}$, $f_{5j}=t^{w_{5j}}$, $f_{6j}=t^{w_{6j}}$.
The Cramer's rule for solving linear systems shows $f_{1j},\ldots,f_{6j}$ to satisfy the conditions~(\ref{eqIIId108}).

\textbf{Case B.} Assume that $j\in\{1,\ldots,q_1\}$, and the minimum over the set $\{w_{4j},w_{5j},w_{6j}\}$ is attained at least twice.
From item 3 it now follows that $\min\{w_{4j},w_{5j},w_{6j}\}=h\leq a.$ Equation~(\ref{eqIIId1}) implies that $h>0$.
We assume w.l.o.g. that $w_{4j}=w_{5j}=h$, $w_{6j}=u\geq h$.
Then the degree of \begin{equation}\label{f5j3}f_{5j}=-\frac{\gamma t^hL_{11}}{L_{12}}-\frac{t^uL_{13}}{L_{12}}+t^{a}\end{equation}
equals $h$ for some $\gamma\in\C^*$. We also set $f_{4j}=\gamma t^h$, $f_{6j}=t^u$.

Now for $i'\in\{1,2,3\}$ we set $$f_{i'j}=\frac{f_{4j}\left(L_{i'1}-\frac{L_{i'2}L_{11}}{L_{12}}\right)+
f_{6j}\left(L_{i'3}-\frac{L_{i'2}L_{13}}{L_{12}}\right)+t^{a}L_{i'2}}{\det L}.$$
Item 5 implies that $f_{1j}=f_{2j}=f_{3j}=0$. From~(\ref{f5j3}) it follows that
$$f_{i'j}=\frac{f_{4j} L_{i'1}+f_{5j}L_{i'2}+f_{6j}L_{i'3}}{\det L},$$
so the Cramer's rule implies the conditions~(\ref{eqIIId108}).

\textbf{Case C.} Now let $j\in\{q_1+1,\ldots,q_1+\ldots+q_4\}$. We will consider only the case when $j\in\{q_1+1,\ldots,q_1+q_2\}$ because
the cases $j\in\{q_1+q_2+1,\ldots,q_1+q_2+q_3\}$ and $j\in\{q_1+q_2+q_3+1,\ldots,q_1+\ldots+q_4\}$ can be considered in the same fashion.
Then by item 2, the tuple $(\deg\lambda_{13},\deg\lambda_{23},\deg\lambda_{33},\deg\lambda_{63})=(b,0,0,0)$ realizes
the tropical dependence of rows of $W[1,2,3,6]$. Now by Lemma~\ref{sys1e}, there exist elements $f_{1j},f_{2j},f_{3j},f_{6j}\in\Kf$ of degrees
$w_{1j},w_{2j},w_{3j},w_{6j}$, respectively, such that $\lambda_{13}f_{1j}+\lambda_{23}f_{2j}+\lambda_{33}f_{3j}+\lambda_{63}f_{6j}=0$.
Items 1 and 3 imply that the elements $f_{4j}=\lambda_{11}f_{1j}+\lambda_{21}f_{2j}+\lambda_{31}f_{3j}$ and
$f_{5j}=\lambda_{12}f_{1j}+\lambda_{22}f_{2j}+\lambda_{32}f_{3j}$ both have a zero degree.
These settings also satisfy the conditions~(\ref{eqIIId108}).

\textbf{Case D.} Let $j\in\{q_1+\ldots+q_4+1,\ldots,n\}$. By the assumptions of case (iii) of Theorem~\ref{hugecases0}, $\Su_j(\P(W))$
has a four-element subset $\{u,v,y,z\}$ that includes no support of a column of $\P(W)$.
Item 5 shows that $\deg\det\Lambda[p,q,r]=0$ for any distinct $p,q,r\in\{u,v,y,z\}$,
items 2 and 4 that $\sum_{h=1}^3\min_{s=1}^6\{\deg \lambda_{sr}+w_{sj}\}=0$.
Thus by Lemma~\ref{sys}, there exist elements $f_{1j},\ldots,f_{6j}\in\Kf$ of the degrees $w_{1j},\ldots,w_{6j}$, respectively,
such that the conditions~(\ref{eqIIId108}) are satisfied.

Cases A--D cover all the possibilities, so there exists a lift $F$ of $W$ such that
the condition~(\ref{eqIIId108}) is satisfied. Thus $rank F\leq3$, so by Definition~\ref{Kapr},
$\Kap(W)\leq3$. The contradiction with item 1 shows that no such $W$ exists.
\end{proof}

\subsection{Case (i)}

To finalize the proof of the main result of this section, we need to consider the case (i) of Theorem~\ref{hugecases0}.
We need to deal with matrices $V\in\R^{6\times n}$ of a more general form, namely, when

\begin{equation}\label{eqIpre1}\P(V)=\left(\begin{array}{c|c|c}
0\ldots0&P_1&\i\ldots\i\\
0\ldots0&&\i\ldots\i\\\hline
\i\ldots\i&0\ldots0&\i\ldots\i\\
\i\ldots\i&0\ldots0&\i\ldots\i\\\hline
\i\ldots\i&P_3&0\ldots0\\
\underbrace{\i\ldots\i}_{u_1}&\underbrace{}_{u_2}&\underbrace{0\ldots0}_{u_3}\end{array}\right),\end{equation}
where $P_1$ and $P_3$ are matrices over $\B$.

We introduce some notation to be used throughout this subsection.
By $L^1$ we denote the set of all tuples $\{l^1=(l_1^1,l_2^1,l_3^1,l_5^1)\}$, $\min\{l_1^1,l_2^1,l_3^1,l_5^1\}=0$,
that realize the tropical dependence of rows of $V[1,2,3,5]$.
By $L^2$ we denote the set of all tuples $\{l^2=(l_1^2,l_3^2,l_4^2,l_5^2)\}$, $\min\{l_1^2,l_3^2,l_4^2,l_5^2\}=0$,
that realize the tropical dependence of rows of $V[1,3,4,5]$.
By $L^3$ we denote the set of all tuples $\{l^3=(l_1^3,l_3^3,l_5^3,l_6^3)\}$, $\min\{l_1^3,l_3^3,l_5^3,l_6^3\}=0$,
that realize the tropical dependence of rows of $V[1,3,5,6]$.

Let $j\in\{1,\ldots,n\}$, $l\in\T^4$. By $\Theta_1(l,j)$ we denote the set of all $\theta_1\in\{1,2,3,5\}$ that provide
the minimum for $\min_{\theta_1}\{\deg v_{\theta_1 j}+l_{\theta_1}\}$.
By $\Theta_2(l,j)$ we denote the set of all $\theta_2\in\{1,3,4,5\}$ that provide
the minimum for $\min_{\theta_2}\{\deg v_{\theta_2 j}+l_{\theta_2}\}$.
By $\Theta_3(l,j)$ we denote the set of all $\theta_3\in\{1,3,5,6\}$ that provide
the minimum for $\min_{\theta_3}\{\deg v_{\theta_3 j}+l_{\theta_3}\}$.

We need the following lemmas.

\begin{lem}\label{cIpl1}
Let a matrix $V\in\R^{6\times n}$ be such that~(\ref{eqIpre1}) holds, and $\trop(V)=3$.
Assume that for every $l^1\in L^1$, $l^3\in L^3$, $j''\in\{u_1+1,\ldots,u_1+u_2\}$ it holds that
$\left|\Theta_1(l^1,j'')\cup\Theta_3(l^3,j'')\right|\geq3$. Assume also that for every $l^2\in L^2$
there exists $l^3\in L^3$ such that for every $j'\in\{1,\ldots,u_1\}$ it holds that
$\left|\Theta_2(l^2,j')\cup\Theta_3(l^3,j')\right|\geq3$. Then $\Kap(V)=3$.
\end{lem}

\begin{proof}
1. From Lemma~\ref{r->b} it follows that $l_1^1=l_2^1=0$  for every $l^1\in L_1$, 
$l_3^2=l_4^2=0$, $l_1^2>0$, $l_5^2>0$ for every $l^2\in L_2$,
$l_5^3=l_6^3=0$ for every $l^3\in L_3$.

2. Lemma~\ref{monot} implies that $\trop(V[1,2,3,4,5])\leq3$. Theorem~\ref{cjrthr} then shows that
$\Kap(V[1,2,3,4,5])\leq3$. By Definition~\ref{Kapr}, there exists $F'\in\Kf^{5\times n}$ such that
$V[1,2,3,4,5]=\deg F'$ and $rank(F')\leq3$.

3. Every four rows of $F'$ are therefore linearly dependent over $\Kf$. In particular, for every $j\in\{1,\ldots,n\}$ it holds that
\begin{equation}\label{Ip1}\lambda_{11} f'_{1j}+\lambda_{21} f'_{2j}+\lambda_{31} f'_{3j}+\lambda_{51} f'_{5j}=0,
\lambda_{12} f'_{1j}+\lambda_{32} f'_{3j}+\lambda_{42} f'_{4j}+\lambda_{52} f'_{5j}=0.\end{equation}
We set also $\lambda_{41}=\lambda_{61}=\lambda_{22}=\lambda_{62}=0$.
Multiplying the equations~(\ref{Ip1}) by elements from $\Kf^*$, we assume w.l.o.g. that
$\min_{\theta=1}^6\{\deg\lambda_{\theta1}\}=\min_{\theta=1}^6\{\deg\lambda_{\theta2}\}=0$.
Lemma~\ref{li->tr} now shows that $l^1=(\deg\lambda_{11} \deg\lambda_{21} \deg\lambda_{31} \deg\lambda_{51})\in L^1$,
$l^2=(\deg\lambda_{12} \deg\lambda_{32} \deg\lambda_{42} \deg\lambda_{52})\in L^2$.

4. By the assumption of the lemma, there exists $l^3\in L^3$ such that $\left|\Theta_2(l^2,j')\cup\Theta_3(l^3,j')\right|\geq3$
for every $j'\in\{1,\ldots,u_1\}$.
Now we set $\lambda_{23}=\lambda_{43}=0$, $\lambda_{13}=\xi_1t^{l^3_{1}}$, $\lambda_{33}=\xi_3t^{l^3_{3}}$,
$\lambda_{53}=\xi_5t^{l^3_{5}}$, $\lambda_{63}=\xi_6t^{l^3_{6}}$.
The infiniteness of $\C$ allows us to find $\xi_\iota\in\C^*$ such that
\begin{equation}\label{Ipre1111}deg(\lambda_{pi'}\lambda_{q3}-\lambda_{qi'}\lambda_{p3})=\min\{\deg \lambda_{pi'}+\deg \lambda_{q3},
\deg \lambda_{qi'}+\deg \lambda_{p3}\}\end{equation}
for every distinct $p,q\in\{1,2,3,4,5,6\}$ and every $i'\in\{1,2\}$.

5. We set $f_{k'''j'''}=f'_{k'''j'''}$ for $k'''\in\{1,2,3,4,5\}$,
$j'''\in\{u_1+u_2+1,\ldots,n\}$. By item 2, $\deg f_{k'''j'''}=v_{k'''j'''}$. Set also
\begin{equation}\label{Ipre111}f_{6j'''}=-\frac{\lambda_{13}f_{1j'''}+\ldots+\lambda_{53}f_{5j'''}}{\lambda_{63}}.\end{equation}
Item 1 and the equation~(\ref{eqIpre1}) imply that $\deg\left(\lambda_{53}f_{5j'''}\right)=0$, and that all the other terms
in the numerator have a positive degree. So we get $\deg f_{6j'''}=0=v_{6j'''}$. Moreover, the equation~(\ref{Ipre111}) implies that the condition
\begin{equation}\label{Ip112}\lambda_{1i} f_{1j}+\lambda_{2i} f_{2j}+\lambda_{3i} f_{3j}+
\lambda_{4i} f_{4j}+\lambda_{5i} f_{5j}+\lambda_{6i} f_{6j}=0\end{equation}
holds for $j\in\{u_1+u_2+1,\ldots,n\}$, $i=3$. The equations~(\ref{Ip1}) also show that~(\ref{Ip112})
holds for $j\in\{u_1+u_2+1,\ldots,n\}$, $i\in\{1,2\}$.

6. Item 4 shows that for $j'\in\{1,\ldots,u_1\}$ the matrix
$$\left(\begin{array}{cc}
\lambda_{12}&\lambda_{13}\\
\lambda_{32}&\lambda_{33}\\
\lambda_{42}&\lambda_{43}\\
\lambda_{52}&\lambda_{53}\\
\lambda_{62}&\lambda_{63}
\end{array}\right)$$
and the tuple $(v_{1j'},v_{3j'},v_{4j'},v_{5j'},v_{6j'})$
satisfy the assumptions of Lemma~\ref{sys2}. Thus there exist $f_{1j'},f_{3j'},f_{4j'},f_{5j'},f_{6j'}\in\Kf$ such that
$\deg f_{k'j'}=v_{k'j'}$ for $k'\in\{1,3,4,5,6\}$, and the conditions~(\ref{Ip112})
hold for $j\in\{1,\ldots,u_1\}$, $i\in\{2,3\}$. We set
\begin{equation}\label{Ip113}f_{2j'}=-\frac{\lambda_{11}f_{1j'}+\lambda_{31}f_{2j'}+\ldots+\lambda_{61}f_{6j'}}{\lambda_{21}}.\end{equation}
Item 1 and the equation~(\ref{eqIpre1}) imply that $\deg\left(\lambda_{11}f_{1j'}\right)=0$, and that all the other terms
in the numerator have a positive degree. So we get $\deg f_{2j'}=0=v_{2j'}$. Moreover,~(\ref{Ip113}) shows that~(\ref{Ip112})
holds for $j\in\{1,\ldots,u_1\}$, $i=1$.

7. The assumption of the lemma shows that
$\left|\Theta_1(l^1,j'')\cup \Theta_3(l^3,j'')\right|\geq3$
for every $j''\in\{u_1+1,\ldots,u_1+u_2\}$. The equations~(\ref{Ipre1111}) then imply that the matrix
$$\left(\begin{array}{cc}
\lambda_{11}&\lambda_{13}\\
\lambda_{21}&\lambda_{23}\\
\lambda_{31}&\lambda_{33}\\
\lambda_{51}&\lambda_{53}\\
\lambda_{61}&\lambda_{63}
\end{array}\right)$$
and the tuple $(v_{1j''},v_{2j''},v_{3j''},v_{5j''},v_{6j''})$
satisfy the assumptions of Lemma~\ref{sys2}. Thus there exist $f_{1j''},f_{2j''},f_{3j''},f_{5j''},f_{6j''}\in\Kf$ such that
$\deg f_{k''j''}=v_{k''j''}$ for $k''\in\{1,2,3,5,6\}$, and the conditions~(\ref{Ip112})
hold for $j\in\{u_1+1,\ldots,u_1+u_2\}$, $i\in\{1,3\}$. We set
\begin{equation}\label{Ip114}f_{4j''}=-\frac{\lambda_{12}f_{1j''}+\lambda_{22}f_{2j''}+\lambda_{32}f_{3j''}+\lambda_{52}f_{5j''}+
\lambda_{62}f_{6j''}}{\lambda_{42}}.\end{equation}
By items 1 and 3, $\deg\left(\lambda_{32}f_{3j''}\right)=0$, and all the other terms
in the numerator have a positive degree. So we get $\deg f_{4j''}=0=v_{4j''}$. Moreover, the equations~(\ref{Ip114}) imply
that~(\ref{Ip112}) holds for $j\in\{u_1+1,\ldots,u_1+u_2\}$, $i=2$.

8. Items 5--7 show that there exists a lift $F$ of $W$ such that
the condition~(\ref{Ip112}) is satisfied for any $j\in\{1,\ldots,n\}$, $i\in\{1,2,3\}$.
Thus $rank F\leq3$, so by Definition~\ref{Kapr}, $\Kap(V)\leq3$.
Now Theorem~\ref{dssthr} shows that $\Kap(V)=3$.
\end{proof}

\begin{lem}\label{cIpl2}
Let a matrix $V\in\R^{6\times n}$ be such that~(\ref{eqIpre1}) holds, and $\trop(V)=3$.
Let for some $l^1\in L^1$, $l^2\in L^2$, $l^3\in L^3$ it holds that
$\left|\Theta_2(l^2,j')\cup \Theta_3(l^3,j')\right|\geq3$ for every $j'\in\{1,\ldots,u_1\}$,
$\left|\Theta_1(l^1,j'')\cup \Theta_3(l^3,j'')\right|\geq3$ for every $j''\in\{u_1+1,\ldots,u_1+u_2\}$,
$\left|\Theta_1(l^1,j''')\cup \Theta_2(l^2,j''')\right|\geq3$ for every $j'''\in\{u_1+u_2+1,\ldots,n\}$.
Then $\Kap(V)=3$.
\end{lem}

\begin{proof}
1. From Lemma~\ref{r->b} it follows that $l_1^1=l_2^1=0$,
$l_3^2=l_4^2=0$, $l_1^2>0$, $l_5^2>0$, $l_5^3=l_6^3=0$.

2. For every $k\in\{1,2,3,4,5,6\}$, $i\in\{1,2,3\}$ we set $$\lambda_{ki}=\zeta_{ki} t^{l^i_k}$$
if the value of $l^i_k$ is defined, and $\lambda_{ki}=0$ otherwise.
The infiniteness of $\C$ allows us to find $\zeta_{ki}\in\C^*$ such that
$$\deg(\lambda_{pi'}\lambda_{qi''}-\lambda_{qi'}\lambda_{pi''})=\min\{\deg \lambda_{pi'}+\deg \lambda_{qi''},
\deg \lambda_{qi'}+\deg \lambda_{pi''}\}$$
for every distinct $p,q\in\{1,2,3,4,5,6\}$ and distinct $i',i''\in\{1,2,3\}$.
The matrix $(\lambda_{ki})$ is denoted by $\Lambda\in\Kf^{6\times3}$.

3. Let $j'\in\{1,\ldots,u_1\}$, $j''\in\{u_1+1,\ldots,u_1+u_2\}$, $j'''\in\{u_1+u_2+1,\ldots,n\}$.
The assumptions of the lemma and item 2 show that the matrix $\Lambda[1,2,3,4,5|1,2]$ and the tuple
$(v_{1j'''},v_{2j'''},v_{3j'''},v_{4j'''},v_{5j'''})$ satisfy the assumptions of Lemma~\ref{sys2}.
We also see that the matrix $\Lambda[1,2,3,5,6|1,3]$ and the tuple
$(v_{1j''},v_{2j''},v_{3j''},v_{5j''},v_{6j''})$ satisfy Lemma~\ref{sys2},
the matrix $\Lambda[1,3,4,5,6|2,3]$ and the tuple
$(v_{1j'},v_{3j'},v_{4j'},v_{5j'},v_{6j'})$ satisfy Lemma~\ref{sys2}.
Thus there exist $f_{k'j'}$, $f_{k''j''}$, $f_{k'''j'''}$ for $k'\in\{1,3,4,5,6\}$, $k''\in\{1,2,3,5,6\}$, $k'''\in\{1,2,3,4,5\}$,
such that $\deg f_{k'j'}=v_{k'j'}$, $\deg f_{k''j''}=v_{k''j''}$, $\deg f_{k'''j'''}=v_{k'''j'''}$, and
and the conditions~(\ref{Ip112}) hold for $i\in\{1,2\}$ if $j\in\{u_1+u_2+1,\ldots,n\}$,
for $i\in\{1,3\}$ if $j\in\{u_1+1,\ldots,u_1+u_2\}$, for $i\in\{2,3\}$ if $j\in\{1,\ldots,u_1\}$.

4. Now we define $f_{2j'}$, $f_{4j''}$, $f_{6j'''}$ by~(\ref{Ip113}), (\ref{Ip114}), (\ref{Ipre111}), respectively.
Note that from item 1 and the equation~(\ref{eqIpre1}) it follows that $\deg f_{2j'}=0=v_{2j'}$, $\deg f_{4j''}=0=v_{4j''}$,
$\deg f_{6j'''}=0=v_{6j'''}$.

5. Items 3--4 show that the matrix $F$ constructed is a lift of $V$, and the conditions~(\ref{Ip112})
hold for every $i\in\{1,2,3\}$, $j\in\{1,\ldots,n\}$. Thus $rank F\leq3$, so by Definition~\ref{Kapr},
$\Kap(V)\leq3$. Now Theorem~\ref{dssthr} implies that $\Kap(V)=3$.
\end{proof}

\begin{thr}\label{caseIpre}
Let a matrix $V\in\R^{6\times n}$ be such that $\trop(V)=3$. Let $\P(V)$ be formed by the columns
$$\left(
\begin{array}{c}
0\\
0\\
\i\\
\i\\
\i\\
\i
\end{array}
\right),\,
\left(
\begin{array}{c}
\i\\
\i\\
0\\
0\\
0\\
0
\end{array}
\right),\,
\left(
\begin{array}{c}
\i\\
\i\\
\i\\
\i\\
0\\
0
\end{array}
\right).$$
Then $\Kap(V)=3$.
\end{thr}

\begin{proof}
1. We have up to p.r.c. that
\begin{equation}\label{eqIp1}
\P(V)=\left(\begin{array}{c|c|c}
0\ldots0&\i\ldots\i&\i\ldots\i\\
0\ldots0&\i\ldots\i&\i\ldots\i\\
\i\ldots\i&0\ldots0&\i\ldots\i\\
\i\ldots\i&0\ldots0&\i\ldots\i\\
\i\ldots\i&0\ldots0&0\ldots0\\
\underbrace{\i\ldots\i}_{u_1}&\underbrace{0\ldots0}_{u_2}&\underbrace{0\ldots0}_{u_3}\end{array}\right).
\end{equation}
By Lemma~\ref{multiply}, we can assume w.l.o.g. that
the minimal element of every column of $V$ equals $0$.

2. Lemma~\ref{r->b} shows that $l_1^1=l_2^1=0$, $l_3^1>0$, $l_5^1>0$ for every $l^1\in L_1$,
$l_3^2=l_4^2=0$, $l_1^2>0$, $l_5^2>0$ for every $l^2\in L_2$, $l_5^3=l_6^3=0$, $l_1^3>0$, $l_3^3\geq0$
for every $l^3\in L_3$.

3. Item 2 and the equation~(\ref{eqIp1}) imply that $6\in\Theta_3(l^3,j'')$
for every $l^3\in L^3$, $j''\in\{u_1+1,\ldots,u_1+u_2\}$. Thus we obtain
$\left|\Theta_1(l^1,j'')\cup \Theta_3(l^3,j'')\right|\geq3$ for every $l^1\in L^1$.

Now we consider the two special cases, A and B.

\textbf{Case A.} Let the cardinality of $L^2$ be different from 1. Since $\trop(V)=3$,
$L^2$ is not empty. Thus there are different elements containing in $L_2$.
Remember that by item 2, $l^2=(l^2_1,0,0,l^2_5)$ for every $l^2\in L^2$.
We have now the four cases to be considered.

A1. Assume that the value of $l^2_1$ is independent on $l^2\in L^2$.
Then for every $l^2\in L^2$ and $j\in\{1,\ldots,n\}$ it holds that
$\left|\Theta_2(l^2,j)\setminus\{5\}\right|\geq 2$.
Item 2 thus shows that $\left|\Theta_1(l^1,j)\cup\Theta_2(l^2,j)\right|\geq3$ for every $l^1\in L^1$.
From item 3 it then follows that $V$ up to p.r.c. satisfies the assumptions of Lemma~\ref{cIpl1}, so $\Kap(V)=3$.

A2. Assume that the value of $l^2_5$ is independent on $l^2\in L^2$.
Then for every $l^2\in L^2$ and $j\in\{1,\ldots,n\}$ it holds that
$\left|\Theta_2(l^2,j)\setminus\{1\}\right|\geq 2$.
Item 2 thus shows that $\left|\Theta_2(l^2,j)\cup\Theta_3(l^3,j)\right|\geq3$ for every $l^3\in L^3$.
From item 3 it then follows that $V$ satisfies the assumptions of Lemma~\ref{cIpl1}, so $\Kap(V)=3$.

A3. Now let the value $l^2_1-l^2_5=d_{15}$ be independent of $l^2\in L^2$.
Then by Lemma~\ref{conve}, for every $(l^2_1,0,0,l^2_5)\in L^2$ there exists a small enough $\varepsilon>0$ such that
either $(l^2_1+\varepsilon,0,0,l^2_5+\varepsilon)\in L^2$ or $(l^2_1-\varepsilon,0,0,l^2_5-\varepsilon)\in L^2$.
This implies that for every $j\in\{1,\ldots,n\}$ it holds that either $\{1,5\}\subset\Theta_2(l^2,j)$ or $\{3,4\}\subset\Theta_2(l^2,j)$.
Thus if $d_{15}\geq0$, then item 2 implies that $\left|\Theta_1(l^1,j)\cup \Theta_2(l^2,j)\right|\geq3$ for every $l^1\in L^1$.
Analogously, if $d_{15}\leq0$, then $\left|\Theta_2(l^2,j)\cup\Theta_3(l^3,j)\right|\geq3$ for every $l^3\in L^3$.
From item 3 it now follows that $V$ up to p.r.c. satisfies the assumptions of Lemma~\ref{cIpl1}, so $\Kap(V)=3$.

A4. Finally, let all the values $l^2_1$, $l^2_5$, $l^2_1-l^2_5$ depend on $l^2\in L^2$. Then by Lemma~\ref{conve}, we can assume w.l.o.g.
that there exist $l^2,\widetilde{l^2},\widehat{l^2},\overline{l^2}\in L^2$ such that $\overline{l^2_1}>l^2_1$, $\widehat{l^2_5}>l^2_5$,
$\widetilde{l^2_1}\geq l^2_1$, $\widetilde{l^2_5}\geq l^2_5$, $\widetilde{l^2_5}-\widetilde{l^2_1}\neq l^2_5-l^2_1$.
Then we set $h=\min\{\widetilde{l^2_1}-l^2_1,\widetilde{l^2_5}-l^2_5\}$, $\ell^2=(l^2_1+h,0,0,l^2_5+h)$.
Lemma~\ref{conve} implies that $\ell^2\in L^2$. Note that for every $j\in\{1,\ldots,n\}$ it holds that either
$\left|\Theta_2(\ell^2,j)\right|\geq 3$ or $\Theta_2(\ell^2,j)=\{3,4\}$.
Thus we see that the matrix $V$ and the tuples $l^1$,$\ell^2$,$l^3$ satisfy the assumptions of Lemma~\ref{cIpl1}, so $\Kap(V)=3$.

The cases A1--A4 cover all the possibilities, so under the assumptions of Case A we obtain $\Kap(V)=3$.

\textbf{Case B.} Now we assume that $L^2$ is a singleton $\{l^2\}$, and that for every $l^1\in L_1$ it holds that $l^1_3<l^1_5$.
We will prove that in this case $\Kap(V)=3$.

B1. By Lemma~\ref{monot}, $\trop(V[1,2,3,4,5])\leq3$. Theorem~\ref{cjrthr} implies that
$\Kap(V[1,2,3,4,5])\leq3$. Definition~\ref{Kapr} now shows that there exists $G''\in\Kf^{5\times n}$
such that $V[1,2,3,4,5]=\deg G''$ and $rank(G'')\leq3$. Analogously, there exists $G'$
such that $V[1,3,4,5,6]=\deg G'$ and $rank(G')\leq3$.

B2. Every four rows of $G'$ and of $G''$ are therefore linearly dependent over $\Kf$. In particular, for every $j\in\{1,\ldots,n\}$ it holds that
\begin{equation}\label{Ip150}\lambda_{11} g''_{1j}+\lambda_{21} g''_{2j}+\lambda_{31} g''_{3j}+\lambda_{51} g''_{5j}=0,
\lambda_{12} g''_{1j}+\lambda_{32} g''_{3j}+\lambda_{42} g''_{4j}+\lambda_{52} g''_{5j}=0,\end{equation}
\begin{equation}
\label{Ip151}\widetilde{\lambda}_{12} g'_{1j}+\widetilde{\lambda}_{32} g'_{3j}+\widetilde{\lambda}_{42} g'_{4j}+
\widetilde{\lambda}_{52} g'_{5j}=0,
\lambda_{13} g'_{1j}+\lambda_{33} g'_{3j}+\lambda_{53} g'_{5j}+\lambda_{63} g'_{6j}=0.
\end{equation}

We also set $\lambda_{41}$$=$$\lambda_{61}$$=$$\lambda_{22}$$=$$\lambda_{62}$$=$$\widetilde{\lambda_{22}}$$=$$\widetilde{\lambda_{62}}$$=
$$\lambda_{23}$$=$$\lambda_{43}$$=$$0$.
Multiplying the equations~(\ref{Ip150}) and~(\ref{Ip151}) by elements from $\Kf^*$, we assume w.l.o.g. that
$$\min_{\theta=1}^6\{\deg\lambda_{\theta1}\}=\min_{\theta=1}^6\{\deg\lambda_{\theta2}\}=
\min_{\theta=1}^6\{\deg\widetilde{\lambda}_{\theta2}\}=\min_{\theta=1}^6\{\deg\lambda_{\theta3}\}=0.$$

B3. From Lemma~\ref{li->tr} it now follows that
$(\deg\lambda_{11} \deg\lambda_{21} \deg\lambda_{31} \deg\lambda_{51})\in L^1$,
$(\deg\lambda_{12} \deg\lambda_{32} \deg\lambda_{42} \deg\lambda_{52})\in L^2$,
$(\deg\widetilde{\lambda}_{12} \deg\widetilde{\lambda}_{32} \deg\widetilde{\lambda}_{42} \deg\widetilde{\lambda}_{52})\in L^2$,
$(\deg\lambda_{13} \deg\lambda_{33} \deg\lambda_{53} \deg\lambda_{63})\in L^3$.

B4. The assumption of Case B shows that
$$(\deg\lambda_{12} \deg\lambda_{32} \deg\lambda_{42} \deg\lambda_{52})=(\deg\widetilde{\lambda}_{12} \deg\widetilde{\lambda}_{32}
\deg\widetilde{\lambda}_{42} \deg\widetilde{\lambda}_{52}),$$
so, multiplying the rows of $G''$ by elements of a zero degree from $\Kf$, we assume w.l.o.g. that
$(\lambda_{12} \lambda_{32} \lambda_{42} \lambda_{52})=(\widetilde{\lambda}_{12} \widetilde{\lambda}_{32}
\widetilde{\lambda}_{42} \widetilde{\lambda}_{52}).$

B5. We also set $g_{k''j''}=g''_{k''j''}$ for $k''\in\{1,2,3,4,5\}$, $j''\in\{u_1+u_2+1,\ldots,n\}$, then $\deg g_{k''j''}=v_{k''j''}$.
Set \begin{equation}\label{Ip155}g_{6j''}=-\frac{\lambda_{13}g_{1j''}+\ldots+\lambda_{53}g_{5j''}}{\lambda_{63}}.\end{equation}
Note that $\deg\left(\lambda_{53}g_{5j''}\right)=0$, and by item 1, all the other terms in the numerator have a positive degree.
So we get $\deg g_{6j''}=0=v_{6j''}$. Moreover, the equation~(\ref{Ip155}) implies that the condition
\begin{equation}\label{Ip156}\lambda_{1i}g_{1j}+\ldots+\lambda_{6i}g_{6j}=0\end{equation}
holds for $j\in\{u_1+u_2+1,\ldots,n\}$, $i=3$. The equation~(\ref{Ip150}) shows that~(\ref{Ip156}) holds for
$j\in\{u_1+u_2+1,\ldots,n\}$, $i\in\{1,2\}$.

B6. We now set $g_{k'j'}=g'_{k'j'}$ for $k'\in\{1,3,4,5,6\}$, $j'\in\{1,\ldots,u_1\}$, then $\deg g_{k'j'}=v_{k'j'}$.
Set \begin{equation}\label{Ip160}g_{2j'}=-\frac{\lambda_{11}g_{1j'}+\lambda_{31}g_{3j'}+\ldots+\lambda_{61}g_{6j'}}{\lambda_{21}}.\end{equation}
Note that $\deg\left(\lambda_{11}g_{1j'}\right)=0$, and by item 1, all the other terms in the numerator have a positive degree.
So we get $\deg g_{1j'}=0=v_{1j'}$. Moreover, the equation~(\ref{Ip160}) implies that the condition~(\ref{Ip156})
holds for $j\in\{1,\ldots,u_1\}$, $i=1$. The equation~(\ref{Ip151}) shows that~(\ref{Ip156}) holds for
$j\in\{1,\ldots,u_1\}$, $i\in\{2,3\}$.

B7. Finally, let $\jmath\in\{u_1+1,\ldots,u_1+u_2\}$.

B7.1. By item 2, $\deg\lambda_{53}=0$. The infiniteness of $\C$ therefore allows us to find $g_{3\jmath}$, $g_{5\jmath}$ such that
$\deg g_{3\jmath}=v_{3\jmath}=0$, $\deg g_{5\jmath}=v_{5\jmath}=0$,  $\deg(\lambda_{33}g_{3\jmath}+\lambda_{53}g_{5\jmath})=0.$

B7.2. By item B3, $(\deg\lambda_{11} \deg\lambda_{21} \deg\lambda_{31} \deg\lambda_{51})\in L^1$,
so from item the assumption of Case B it follows that $\deg(\lambda_{31}g_{3\jmath})<\deg(\lambda_{51}g_{5\jmath})$.
The definition of the set $L^1$ now implies that
either $\tau=1$ or $\tau=2$ provides the minimum for $\min\limits_{\tau\in\{1,2,3,5\}}\{\deg\lambda_{\tau1}+v_{\tau\jmath}\}$.
The infiniteness of $\C$ therefore allows us to find $g_{1\jmath}$, $g_{2\jmath}$ such that $\deg g_{1\jmath}=v_{1\jmath}$,
$\deg g_{2\jmath}=v_{2\jmath}$, $\lambda_{11}g_{1\jmath}+\lambda_{21}g_{2\jmath}+\lambda_{31}g_{3\jmath}+\lambda_{51}g_{5\jmath}=0$.
By item B2, $\lambda_{41}=\lambda_{61}=0$, so the condition~(\ref{Ip156})
holds for $j\in\{u_1+1,\ldots,u_1+u_2\}$, $i=1$ for any $g_{4\jmath}$, $g_{6\jmath}$.

B7.3. We set
\begin{equation}\label{Ip161}g_{6\jmath}=-\frac{\lambda_{13}g_{1\jmath}+
\lambda_{33}g_{3\jmath}+\lambda_{53}g_{5\jmath}}{\lambda_{63}}.\end{equation}
By item 1, $\deg\left(\lambda_{13}g_{1\jmath}\right)>0$, so from item B7.1 it follows that $\deg g_{6\jmath}=0=v_{6\jmath}$.
By item B2, $\lambda_{23}=\lambda_{43}=0$, so the condition~(\ref{Ip156})
holds for $j\in\{u_1+1,\ldots,u_1+u_2\}$, $i=3$ for any $g_{4\jmath}$.

B7.4. We set
\begin{equation}\label{Ip162}g_{4\jmath}=-\frac{\lambda_{12}g_{1\jmath}+\lambda_{22}g_{2\jmath}+
\lambda_{32}g_{3\jmath}+\lambda_{52}g_{5\jmath}+\lambda_{62}g_{6\jmath}}{\lambda_{42}}.\end{equation}
Note that $\deg\left(\lambda_{32}g_{3\jmath}\right)=0$, and item 2 shows that all the other terms
in the numerator have a positive degree. So we obtain $\deg g_{4\jmath}=0=v_{4\jmath}$,
and the equation~(\ref{Ip162}) shows that~(\ref{Ip156}) holds for $j\in\{u_1+1,\ldots,u_1+u_2\}$, $i=2$.

B8. Items B5--B7 construct a lift $G$ of $V$ such that the conditions~(\ref{Ip156})
hold for every $i\in\{1,2,3\}$, $j\in\{1,\ldots,n\}$. Thus $rank G\leq3$, so by Definition~\ref{Kapr},
$\Kap(V)\leq3$. Now Theorem~\ref{dssthr} implies that $\Kap(V)=3$ and completes the consideration of Case B.

Now we can \textbf{finalize} the proof of Theorem~\ref{caseIpre}. This is done by \textit{reductio ad absurdum}.

F0. Indeed, assume $\Kap(V)\neq3$.

F1. Cases A and B show that $L^2=\{l^2\}$, and that $\ell^1_5\leq\ell^1_3$ for some $\ell^1\in L^1$.

F2. If for every $\widetilde{j}\in\{u_1+u_2+1,\ldots,n\}$
it holds that $\left|\Theta_1(\ell^1,\widetilde{j})\cup\Theta_2(l^2,\widetilde{j})\right|\geq3$,
then item 3 shows that $V$ satisfies up to p.r.c. the assumptions of Lemma~\ref{cIpl1}.
This gives a contradiction with assumption F0.

Thus there exists $\widetilde{j}\in\{u_1+u_2+1,\ldots,n\}$
such that $\left|\Theta_1(\ell^1,\widetilde{j})\cup\Theta_2(l^2,\widetilde{j})\right|=2$. Items 2 and F1 now imply that
$l^2=(b,0,0,a+b)$ for some $a,b>0$, and $\Theta_1(\ell^1,\widetilde{j})=\Theta_2(l^2,\widetilde{j})=\{1,5\}$.
The $\widetilde{j}$th column of $V$ in this case has the form
$$\left(
\begin{array}{c}
a\\
\alpha_1\\
\gamma_1\\
\gamma_2\\
0\\
0
\end{array}
\right),
$$
where $\alpha_1>a$, $\gamma_1>a+b$, $\gamma_2>a+b$.

F3. Analogously, for $l^3\in L^3$ there exists $\widetilde{j}'\in\{1,\ldots,u_1\}$
such that $\left|\Theta_2(l^2,\widetilde{j}')\cup\Theta_3(l^3,\widetilde{j}')\right|=2$.
Item 2 now implies that $\Theta_2(l^2,\widetilde{j}')=\Theta_3(l^3,\widetilde{j}')=\{1,3\}$,
so in this case the $\widetilde{j}'$th column of $V$ has the form
$$\left(
\begin{array}{c}
0\\
0\\
b\\
\beta_1\\
\nu_1\\
\nu_2
\end{array}
\right),
$$
where $\beta_1>b$, $\nu_1>b$, $\nu_2>b$.

F4. All the matrices that can be obtained from $V$ by permutations of the first two and of the second two rows
also satisfy the assumptions of the theorem being proved.

F5. Item F4 allows us to apply the result of items F2--F3 to the matrices considered in item F4,
thus we see that the following columns appear in $V$:

\begin{equation}\label{123123123}\left(
\begin{array}{c}
0\\
0\\
b\\
\beta_1\\
\nu_1\\
\nu_2
\end{array}
\right),
\left(
\begin{array}{c}
0\\
0\\
\beta_2\\
b\\
\nu_3\\
\nu_4
\end{array}
\right),
\left(
\begin{array}{c}
0\\
0\\
b'\\
\beta'_1\\
\nu_5\\
\nu_6
\end{array}
\right),
\left(
\begin{array}{c}
0\\
0\\
\beta'_2\\
b'\\
\nu_7\\
\nu_8
\end{array}
\right),
\left(
\begin{array}{c}
a\\
\alpha\\
\gamma_1\\
\gamma_2\\
0\\
0
\end{array}
\right),
\left(
\begin{array}{c}
\alpha'\\
a'\\
\gamma'_1\\
\gamma'_2\\
0\\
0
\end{array}
\right)
\end{equation}
where $a'>0$, $b'>0$, $\alpha>a$, $\alpha'>a'$, $\gamma_1>a+b$, $\gamma_2>a+b$, $\gamma'_1>a'+b'$, $\gamma'_2>a'+b'$,
$\beta_1>b$, $\beta_2>b$, $\beta'_1>b'$, $\beta'_2>b'$, $\nu_1,\ldots,\nu_4$ are greater than $b$, $\nu_5,\ldots,\nu_8$ are greater than $b'$.

F6. Item F4 allows us to assume w.l.o.g. that $a\leq a'$, $b\leq b'$. We consider the matrix
formed by the $1$st, $3$rd, $4$th, and $5$th rows of the $1$st, $2$nd, $5$th, and $6$th columns of~(\ref{123123123}):
$$S=\left(
\begin{array}{cccc}
0&0&a&\alpha'\\
b&\beta_2&\gamma_1&\gamma'_1\\
\beta_1&b&\gamma_2&\gamma'_2\\
\nu_1&\nu_3&0&0
\end{array}
\right).
$$

We note that the permanent of $S$ equals $a+2b$, and the minimum in~(\ref{permdef}) is provided
by the unique permutation $(132)\in\S_4$. By Definition~\ref{troprk}, $\trop(V)\geq4$.
The contradiction obtained shows that the assumption F0 fails to hold. Thus we actually have $\Kap(V)=3$.
\end{proof}

\begin{thr}\label{caseI}
Case (i) of Theorem~\ref{hugecases0} is not realizable.
\end{thr}

\begin{proof}
1. Let a matrix $V$ realize case (i), then $\trop(V)=3$, $\Kap(V)>3$, and
\begin{equation}\label{eqI1}\P(V)=\left(\begin{array}{c|c|c}
0\ldots0&\i\ldots\i&\i\ldots\i\\
0\ldots0&\i\ldots\i&\i\ldots\i\\
\i\ldots\i&0\ldots0&\i\ldots\i\\
\i\ldots\i&0\ldots0&\i\ldots\i\\
\i\ldots\i&\i\ldots\i&0\ldots0\\
\underbrace{\i\ldots\i}_{u_1}&\underbrace{\i\ldots\i}_{u_2}&\underbrace{0\ldots0}_{u_3}\end{array}\right).\end{equation}
By Lemma~\ref{multiply}, we can assume w.l.o.g. that the minimal element of every column of $V$ equals $0$.

2. For $i\in\{1,2,3\}$ we set
\begin{equation}\label{eqI2}\mu_i=\min_{j\in\{1,\ldots,n\}: v_{2i-1,j}\neq v_{2i,j}}\left\{\min\left\{v_{2i-1,j},v_{2i,j}\right\}\right\}.\end{equation}
Note that if $\mu_i=\i$ $\left(\mbox{i.e., if $v_{2i-1,j}=v_{2i,j}$ for any $j\in\{1,\ldots,n\}$}\right)$, then the
$(2i-1)$-th and $2i$-th rows of $V$ coincide. Then Theorem~\ref{cjrthr} implies that
$\Kap(V)=\trop(V)$ and gives a contradiction. Thus $\mu_i\in\R$, item 1 implies that $\mu_i>0$.

3. By $g(i)$ we denote any $j\in\{1,\ldots,n\}$ that provides the minimum for~(\ref{eqI2}),
i.e. $\mu_i=\min\{v_{2i-1,g(i)},v_{2i,g(i)}\},$ $v_{2i-1,g(i)}\neq v_{2i,g(i)}$.

4. After renumbering rows and columns of $V$ we can assume w.l.o.g. that
$g(1)\in\{u_1+1,\ldots,u_1+u_2\}$, $g(3)\in\{1,\ldots,u_1\}$.

Assume that $\min\{v_{5j'},v_{6j'}\}<\mu_3$ for every
$j'\in\{u_1+1,\ldots,u_1+u_2\}$. Then by item 2,
$v_{5j'}=v_{6j'}<\mu_3$ for every $j'\in\{u_1+1,\ldots,u_1+u_2\}$,
so we have $$\min_{j'=u_1+1}^{u_1+u_2}\{v_{5j}\}=\overline{\mu_3}<\mu_3.$$
Now we add $-\overline{\mu_3}$ to every element of the last two rows of $V$
and a small enough $-\varepsilon<0$ to every element of the first two.
The matrix obtained then satisfies the assumptions of Theorem~\ref{caseIpre} and gives a contradiction with item 1.

Therefore we have that $\min\{v_{5h(3)},v_{6h(3)}\}\geq\mu_3$ for some $h(3)\in\{u_1+1,\ldots,u_1+u_2\}$.
Analogously, we see that $\min\{v_{1h(1)},v_{2h(1)}\}\geq\mu_1$ for some $h(1)\in\{u_1+u_2+1,\ldots,n\}$.

5. From item 1 it follows that $v_{3g}=v_{4g}=0$ for $g\in\{u_1+1,\ldots,u_1+u_2\}$.
By item 4, $v_{3g(1)}=v_{4g(1)}=0$, and the numbers $v_{5g(1)}$ and
$v_{6g(1)}$ are positive. Analogously, $v_{1g(3)}=v_{2g(3)}=v_{5h(1)}=v_{6h(1)}=0$,
$v_{3g(3)}>0$, $v_{4g(3)}>0$.

6. Assume that $V$ contains columns (we denote the numbers of these columns by $\gamma_1$ and $\gamma_2$, respectively)
of both of the following forms:
\begin{equation}\label{eqI4}\left(\begin{array}{c}
\mu_1\\
\psi_1\\
0\\
0\\
\omega_1\\
\omega_2
\end{array}\right),\end{equation}
\begin{equation}\label{eqI5}\left(\begin{array}{c}
\psi_2\\
\mu_1\\
0\\
0\\
\omega_3\\
\omega_4
\end{array}\right),\end{equation}
where $\psi_1>\mu_1$, $\psi_2>\mu_1$, and $\omega_1,\ldots,\omega_4$ are all greater than $\mu_1+\mu_3$.
Let us note that by item 5, the matrix
$V[5,6,1,2|n,g(3),\gamma_1,\gamma_2]$ equals
$$S=\left(\begin{array}{cccc}
0&m'&\omega_1&\omega_3\\
0&m''&\omega_2&\omega_4\\
\chi_1&0&\mu_1&\psi_2\\
\chi_2&0&\psi_1&\mu_1
\end{array}\right),$$
where $\chi_1$, $\chi_2$ are positive. From item 3 it also follows that $m'\neq m''$ and $\min\{m',m''\}=\mu_3$.
The permanent of $S$ equals $2\mu_1+\mu_3$, and the minimum in~(\ref{permanent}) is given by a unique permutation.
Thus $S$ is tropically non-singular, $\trop(V)\geq4$. The contradiction shows that at most one
of~(\ref{eqI4}) and~(\ref{eqI5}) appears as a column of $V$. We assume w.l.o.g. that~(\ref{eqI4}) does not appear.

7. Let the tuple $(l_{1},l_{2},l_{3},l_{5})$, $\min\{l_{1},l_{2},l_{3},l_{5}\}=0$,
realize the tropical dependence of rows of $V[1,2,3,5]$.
Lemma~\ref{r->b} implies that $l_{1}$$=$$l_{2}$$=$$0$, $l_{3}>0$, $l_{5}>0$.

Assume that $l_{3}\neq\mu_1$. Since the minimum over
$\{v_{1g(1)}+l_{1},v_{2g(1)}+l_{2},v_{3 g(1)}+l_{3},v_{5 g(1)}+l_{5}\}=\{v_{1g(1)},v_{2g(1)},l_{3},v_{5 g(1)}+l_{5}\}$
is attained at least twice, item 3 shows that $l_{5}<\min\{l_{3},\mu_1\}$.
Now from item 5 it follows that $v_{5 h(1)}+l_{5}<\min\{l_{3},\mu_1\}$, from item 4 that
$\min\{v_{1h(1)}+l_{1}, v_{2h(1)}+l_{2}, v_{3 h(1)}+l_{3}\}\geq\min\{l_{3},\mu_1\}$.

The contradiction with Definition~\ref{tropdep} shows that $l_{3}=\mu_1$.
If $l_5<\mu_1$, then item 5 implies that $v_{5 h(1)}+l_5<\mu_1$, item 4 that
$\min\{v_{1h(1)}, v_{2h(1)}, v_{3 h(1)}+\mu_1\}\geq\mu_1$.
This also contradicts Definition~\ref{tropdep}, so $l_5\geq\mu_1$.

8. From Theorem~\ref{tropthr} it follows that every four rows of $V$ are tropically dependent.
Item 7 shows that there exists $y_1\in\R$, $y_1\geq\mu_1$, such that the tuple $(0,0,\mu_1,y_1)$ realizes the tropical
dependence of rows of $V[1,2,3,5]$. The similar argument shows that there exists $y_3\in\R$, $y_3\geq\mu_3$, such that
the tuple $(\mu_3,y_3,0,0)$ realizes the tropical dependence of rows of $V[1,3,5,6]$.

9. There are the two possible cases for the value of $g(2)$.
Case A, $g(2)\in\{u_1+u_2+1,\ldots,n\}$, and Case B, $g(2)\in\{1,\ldots,u_1\}$.
We will consider these cases separately.

{\bf Case A.} So, let $g(2)\in\{u_1+u_2+1,\ldots,n\}$.

A1. The argument similar to one of item 4 shows that
$\min\{v_{3h(2)},v_{4h(2)}\}\geq\mu_2$ for some $h(2)\in\{1,\ldots,u_1\}$.

A2. The argument similar to one of items 7--8 now shows that there exists $y_2\in\R$,
$y_2\geq\mu_2$, such that the tuple $(y_2,0,0,\mu_2)$ realizes the tropical
dependence of rows of $V[1,3,4,5]$.

A3. The argument similar to one of item 6 allows us to assume
w.l.o.g. that the following columns do not appear in $V$:
$$\left(\begin{array}{c}
\mu_1\\
\psi'\\
0\\
0\\
\omega'_1\\
\omega'_2
\end{array}\right),
\left(\begin{array}{c}
0\\
0\\
\omega''_1\\
\omega''_2\\
\mu_3\\
\psi''
\end{array}\right),
\left(\begin{array}{c}
\omega'''_1\\
\omega'''_2\\
\mu_2\\
\psi'''\\
0\\
0
\end{array}\right),
$$
where $\psi'>\mu_1$, $\psi''>\mu_3$, $\psi'''>\mu_2$, $\omega'_1>\mu_1+\mu_3$, $\omega'_2>\mu_1+\mu_3$,
$\omega''_1>\mu_2+\mu_3$, $\omega''_2>\mu_2+\mu_3$, $\omega'''_1>\mu_1+\mu_2$, $\omega'''_2>\mu_1+\mu_2$.

A4. We set
$$\Lambda=\left(\begin{array}{ccc}
1&t^{y_2}&2t^{\mu_3}\\
1&0&0\\
2t^{\mu_1}&1&t^{y_3}\\
0&1&0\\
t^{y_1}&2t^{\mu_2}&1\\
0&0&1
\end{array}\right).$$

We will show that there exists a lift $F$ of $V$ such that
\begin{equation}\label{eqI6}\lambda_{1i}f_{1j}+\lambda_{2i}f_{2j}+\lambda_{3i}f_{3j}+\lambda_{4i}f_{4j}+\lambda_{5i}f_{5j}+\lambda_{6i}f_{6j}=0\end{equation}
for every $i\in\{1,2,3\}$, $j\in\{1,\ldots,n\}$.

A5. First, let $j\in\{1,\ldots,u_1\}$. By item 1, $v_{1j}=v_{2j}=0$, item A3 shows that $v_{5j}\neq\mu_3$, or $v_{6j}\leq\mu_3$,
or $v_{3j}\leq\mu_2+\mu_3$, or $v_{4j}\leq\mu_2+\mu_3$. Items 8 and A2 imply that in any of these cases the matrix
$$\left(\begin{array}{ccc}
t^{y_2}&2t^{\mu_3}\\
1&t^{y_3}\\
1&0\\
2t^{\mu_2}&1\\
0&1
\end{array}\right)$$
(which is obtained from $\Lambda$ by removing the first column and the second row) and the tuple $(v_{1j},v_{3j},v_{4j},v_{5j},v_{6j})$
satisfy the assumptions of Lemma~\ref{sys2}. Thus there exist $f_{1j},f_{3j},f_{4j},f_{5j},f_{6j}\in\Kf$ such that
$\deg f_{k'j}=v_{k'j}$ for $k'\in\{1,3,4,5,6\}$, and
$\lambda_{1i'}f_{1j}+\lambda_{3i'}f_{3j}+\lambda_{4i'}f_{4j}+\lambda_{5i'}f_{5j}+\lambda_{6i'}f_{6j}=0$ for $i'\in\{2,3\}$.
Since $\lambda_{22}=\lambda_{23}=0$, the conditions~(\ref{eqI6}) hold for $i\in\{2,3\}$ and for any $f_{2j}$.

We set \begin{equation}\label{caseI40811}f_{2j}=-\lambda_{11}f_{1j}-\lambda_{31}f_{3j}-\lambda_{41}f_{4j}-\lambda_{51}f_{5j}-\lambda_{61}f_{6j}.\end{equation}
Now the conditions~(\ref{eqI6}) also hold for $i=1$. The equation~(\ref{eqI1}) shows that $\lambda_{11}f_{1j}$ is the unique term with non-positive degree
in the right-hand side of~(\ref{caseI40811}). Thus we have $\deg f_{2j}=0=v_{2j}$.

A6. We note that the cases $j\in\{1,\ldots,u_1\}$, $j\in\{u_1+1,\ldots,u_1+u_2\}$, and $j\in\{u_1+u_2+1,\ldots,n\}$ are
the same up to renumbering the rows and the columns of $V$. Thus for every $j\in\{1,\ldots,n\}$ we can prove that
there exist $f_{1j},\ldots,f_{6j}\in\Kf$ with degrees $v_{1j},\ldots,v_{6j}$, respectively, such that
the conditions~(\ref{eqI6}) hold for every $i\in\{1,2,3\}$.

The definition of $\Lambda$ shows that every row of the matrix $F$ constructed is a linear combination of its first, third, and fifth rows.
By Definition~\ref{Kapr}, $\Kap(V)\leq3$. The contradiction with item 1 completes the consideration of Case A.

{\bf Case B.} Now let $g(2)\in\{1,\ldots,u_1\}$.

B1. The argument similar to one of item 4 shows that
$\min\{v_{3h(2)},v_{4h(2)}\}\geq\mu_2$ for some $h(2)\in\{u_1+u_2+1,\ldots,n\}$.

B2. By Lemma~\ref{monot}, $\trop(V[1,2,3,4,5])\leq3$. Theorem~\ref{cjrthr} implies that
$\Kap(V[1,2,3,4,5])\leq3$. Definition~\ref{Kapr} now shows that there exists $F'\in\Kf^{5\times n}$
such that $V[1,2,3,4,5]=\deg F'$ and $rank(F')\leq3$.

Every four rows of $F'$ are therefore linearly dependent over $\Kf$. In particular, for every $j\in\{1,\ldots,n\}$ it holds that
\begin{equation}\label{IB1}\lambda_{11} f'_{1j}+\lambda_{21} f'_{2j}+\lambda_{31} f'_{3j}+\lambda_{51} f'_{5j}=0,
\lambda_{12} f'_{1j}+\lambda_{32} f'_{3j}+\lambda_{42} f'_{4j}+\lambda_{52} f'_{5j}=0.\end{equation}

B3. Multiplying the equations~(\ref{IB1}) by elements from $\Kf^*$, we assume w.l.o.g.
that $\min\{\deg\lambda_{11},\deg\lambda_{21},\deg\lambda_{31},\deg\lambda_{51}\}=0$,
$\min\{\deg\lambda_{12},\deg\lambda_{32},\deg\lambda_{42},\deg\lambda_{52}\}=0$.
Lemma~\ref{li->tr} shows that the tuples $(\deg\lambda_{11} \deg \lambda_{21} \deg \lambda_{31} \deg \lambda_{51})$ and
$(\deg\lambda_{12} \deg \lambda_{32} \deg \lambda_{42} \deg \lambda_{52})$ realize the tropical dependence of rows of
$V[1,2,3,5]$ and $V[1,3,4,5]$, respectively. Item 7 then implies that $\deg\lambda_{11}=\deg\lambda_{21}=0$,
$\deg\lambda_{31}=\mu_1$, $\deg\lambda_{51}=y_1\geq\mu_1$. Analogously, we can get that
$\deg\lambda_{32}=\deg\lambda_{42}=0$, $\deg\lambda_{12}=\mu_2$, $\deg\lambda_{52}=y_2\geq\mu_2$.
We set also $\lambda_{41}=\lambda_{61}=\lambda_{22}=\lambda_{62}=0$.

B4. Item 8 shows that the tuple $(\mu_3,y_3,0,0)$ realizes the tropical dependence of rows of $V[1,3,5,6]$.
We set $\lambda_{23}=\lambda_{43}=0$. Set also
$\lambda_{13}=\xi_1t^{\mu_3}$, $\lambda_{33}=\xi_3t^{y_3}$, $\lambda_{53}=\xi_5$, $\lambda_{63}=\xi_6$,
where $\xi_\iota\in\C^*$ are such that
\begin{equation}\label{IB11}\deg(\lambda_{p1}\lambda_{q3}-\lambda_{q1}\lambda_{p3})=\min\{\deg \lambda_{p1}+\deg \lambda_{q3},
\deg \lambda_{q1}+\deg \lambda_{p3}\}\end{equation} for every distinct $p,q\in\{1,2,3,4,5,6\}$.
The existence of such $\{\xi_\iota\}$ follows from the infiniteness of $\C$.

B5. Let $j'\in\{u_1+u_2+1,\ldots,n\}$. We set $f_{k'j'}=f'_{k'j'}$ for $k'\in\{1,2,3,4,5\}$.
Item B2 shows that $\deg f_{k'j'}=v_{k'j'}$. Set also
\begin{equation}\label{IB2}f_{6j'}=-\frac{\lambda_{13}f_{1j'}+\lambda_{23}f_{2j'}+\ldots+\lambda_{53}f_{5j'}}{\lambda_{63}}.\end{equation}
Note that $\deg\left(\lambda_{53}f_{5j'}\right)=0$, and, by the equation~(\ref{eqI1}), all the other terms in the numerator have a positive degree.
So we get $\deg f_{6j'}=0=v_{6j'}$. Moreover, the equation~(\ref{IB2}) implies that the condition
\begin{equation}\label{IB3}\lambda_{1i} f_{1j}+\lambda_{2i} f_{2j}+\lambda_{3i} f_{3j}+
\lambda_{4i} f_{4j}+\lambda_{5i} f_{5j}+\lambda_{6i} f_{6j}=0\end{equation}
holds for $j\in\{u_1+u_2+1,\ldots,n\}$, $i=3$. The equations~(\ref{IB1}) show that~(\ref{IB3}) holds for
$j\in\{u_1+u_2+1,\ldots,n\}$, $i\in\{1,2\}$.

B6. Now let $j''\in\{1,\ldots,u_1\}$. The equation~(\ref{IB11}) and items B3--B4 imply that the matrix
$$\left(\begin{array}{cc}
\lambda_{12}&\lambda_{13}\\
\lambda_{32}&\lambda_{33}\\
\lambda_{42}&\lambda_{43}\\
\lambda_{52}&\lambda_{53}\\
\lambda_{62}&\lambda_{63}
\end{array}\right)$$
and the tuple $(v_{1j''},v_{3j''},v_{4j''},v_{5j''},v_{6j''})$ satisfy the conditions of Lemma~\ref{sys2}.
Thus there exist $f_{1j''},f_{3j''},f_{4j''},f_{5j''},f_{6j''}\in\Kf$
such that $\deg f_{k''j''}=v_{k''j''}$ for every $k''\in\{1,3,4,5,6\}$,
and the condition~(\ref{IB3}) holds for $j\in\{1,\ldots,u_1\}$, $i\in\{2,3\}$.
Set also
\begin{equation}\label{IB5}f_{2j''}=-\frac{\lambda_{11}f_{1j''}+\lambda_{31}f_{2j''}+\ldots+\lambda_{61}f_{5j''}}{\lambda_{21}}.\end{equation}
Note that $\deg\left(\lambda_{11}f_{1j''}\right)=0$, and, by the equation~(\ref{eqI1}), all the other terms in the numerator have a positive degree.
So we get $\deg f_{2j''}$$=$$0$$=$$v_{2j''}$. Moreover, the equation~(\ref{IB5}) implies that the condition~(\ref{IB3})
holds also for $j\in\{1,\ldots,u_1\}$, $i=1$.

B7. Finally, let $j'''\in\{u_1+1,\ldots,u_1+u_2\}$. Item 1 implies that $v_{3j'''}=v_{4j'''}=0$,
item 6 that $v_{1j'''}\neq\mu_1$, or $v_{2j'''}\leq\mu_1$,
or $v_{5j'''}\leq\mu_1+\mu_3$, or $v_{6j'''}\leq\mu_1+\mu_3$.
By items B3--B4, in any of these cases the matrix
$$\left(\begin{array}{cc}
\lambda_{11}&\lambda_{13}\\
\lambda_{21}&\lambda_{23}\\
\lambda_{31}&\lambda_{33}\\
\lambda_{51}&\lambda_{53}\\
\lambda_{61}&\lambda_{63}
\end{array}\right)$$
and the tuple $(v_{1j'''},v_{2j'''},v_{3j'''},v_{5j'''},v_{6j'''})$ satisfy the assumptions
of Lemma~\ref{sys2}, so there exist $f_{1j'''},f_{2j'''},f_{3j'''},f_{5j'''},f_{6j'''}$
such that $\deg f_{k'''j'''}=v_{k'''j'''}$ for every $k'''\in\{1,2,3,5,6\}$,
and the condition~(\ref{IB3}) holds for $j\in\{u_1+1,\ldots,u_1+u_2\}$, $i\in\{1,3\}$.
Set also
\begin{equation}\label{IB6}f_{4j'''}=-\frac{\lambda_{12}f_{1j'''}+\lambda_{22}f_{2j'''}+\lambda_{32}f_{3j'''}+
\lambda_{52}f_{5j'''}+\lambda_{62}f_{6j'''}}{\lambda_{42}}.\end{equation}
Note that $\deg\left(\lambda_{32}f_{3j'''}\right)=0$, and, by the equation~(\ref{eqI1}), all the other terms in the numerator have a positive degree.
So we get $\deg f_{4j'''}=0=v_{4j'''}$. Moreover, the equation~(\ref{IB6}) implies that the condition~(\ref{IB3})
holds also for $j\in\{u_1+1,\ldots,u_1+u_2\}$, $i=2$.

B8. Items B5--B7 show the existence of a matrix $F$ such that $V=\deg F$, and the conditions~(\ref{IB3})
hold for every $i\in\{1,2,3\}$, $j\in\{1,\ldots,n\}$. By items B3 and B4, every column of $F$ is then
a linear combination of its first, third, and fifth columns, so by Definition~\ref{Kapr}, $\Kap(V)\leq3$.
The contradiction with item 1 shows that Case B is also not realizable. The proof is complete.
\end{proof}

Now we can prove the main result of this section.

\begin{thr}\label{6x6r3}
Let $A\in\R^{6\times n}$ be such that $\trop(A)=3$. Then $\Kap(A)=3$.
\end{thr}

\begin{proof}
If $\Kap(A)>3$, then, by Theorem~\ref{hugecases0}, we can assume w.l.o.g. that $A$ satisfies one of the cases (i)--(v).
Theorems~\ref{case VI},~\ref{caseII},~\ref{cases IV-V},~\ref{caseIIId},~\ref{caseI}
show that none of these cases is realizable, so $\Kap(A)\leq3$.
By Theorem~\ref{dssthr}, we get $\Kap(A)=3$.
\end{proof}

\section{The main result}
This section finalizes the proof of the main result of our paper.
Example~\ref{6x6} and the following will be now important.

\begin{ex}\cite{DSS}\label{dssex}
Let
$$C=\left(
\begin{array}{ccccccc}
1&1&0&1&0&0&0\\
0&1&1&0&1&0&0\\
0&0&1&1&0&1&0\\
0&0&0&1&1&0&1\\
1&0&0&0&1&1&0\\
0&1&0&0&0&1&1\\
1&0&1&0&0&0&1
\end{array}
\right).$$
Then $\trop(C)=3$, $\Kap(C)=4$.
\end{ex}

\begin{proof}
The proof can be given by a straightforward application of Definitions~\ref{troprk} and~\ref{Kapr}.
See also~\cite[Section 7]{DSS}.
\end{proof}

\begin{thr}\label{extr}
Let $A\in\R^{d\times n}$ be such that $\trop(A)=r$, $\Kap(A)>r$. Then there exist $A'\in\R^{(d+1)\times n}$,
$A''\in\R^{d\times (n+1)}$ such that $\trop(A')=\trop(A'')=r$, $\Kap(A')=\Kap(A'')>r$.
\end{thr}

\begin{proof}
It is enough to note that by Definitions~\ref{troprk} and~\ref{Kapr}, the tropical and Kapranov ranks
of a matrix are invariant with respect to adding a repeating row or column.
\end{proof}

Using the construction provided in~\cite[Section 7]{DSS}, we prove the following theorem.

\begin{thr}\label{extr2}
Let $A\in\R^{d\times n}$ be such that $\trop(A)=r$, $\Kap(A)>r$. Then there exists
$B\in\R^{(d+1)\times(n+1)}$ such that $\trop(B)=r+1$, $\Kap(B)>r+1$.
\end{thr}

\begin{proof}
By Lemma~\ref{multiply}, we can assume w.l.o.g. that the minimal element
of every row and every column of $A$ equals $0$.
By $P$ we denote the largest tropical permanent over all $r$-by-$r$ submatrices of $A$.
We set
$$B=\left(
\begin{array}{ccc|c}
&&&P+1+a_{1n}\\
&A&&\ldots\\
&&&P+1+a_{dn}\\\hline
P+1+a_{d1}&\ldots&P+1+a_{dn}&0\\
\end{array}
\right).$$
Lemma~\ref{multiply} implies that $\trop(B)\leq r+1$. The choice of $P$ also shows that the matrix
$B[h_1,\ldots,h_r,d+1|c_1,\ldots,c_r,n+1]$ is tropically non-singular if a matrix $A[h_1,\ldots,h_r|c_1,\ldots,c_r]$ is.
Thus we see that $\trop(B)=r+1$.

It remains to check that $\Kap(B)>r+1$. Indeed, let $F$ be a lift of $B$.
We denote  $$D=\textbf{I}_{d+1}-\sum\limits_{k=1}^d\frac{f_{k,n+1}\textbf{U}_{k,d+1}}{f_{d+1,n+1}}\in\Kf^{(d+1)\times(d+1)},$$
where $\textbf{I}_{d+1}$ is the identity matrix, $\textbf{U}_{k,d+1}$ the matrix units.
Since $D$ is non-singular, we obtain $rank(DF)=rank(F)$.
One can note that
\begin{equation}\label{metka3sec}
DF=\left(
\begin{array}{ccc|c}
&&&0\\
&\overline{A}&&\ldots\\
&&&0\\\hline
f_{d+1,1}&\ldots&f_{d+1,n}&f_{d+1,n+1}\\
\end{array}
\right),
\end{equation}
where $\overline{A}$ is a lift of $A$. The assumptions of the theorem show that $rank(\overline{A})>r$.
Equation~(\ref{metka3sec}) implies that $rank(DF)=rank(\overline{A})+1$. Thus we see that $rank(F)>r+1$.
\end{proof}

Now we can obtain an answer for Question~\ref{questi}.

\begin{thr}\label{mainthr2}
Let $d,n,r$ be positive integers, $r\leq\min\{d,n\}$. Then $d$-by-$n$ matrices
with tropical rank less than $r$ always have the Kapranov rank less than $r$
if and only if one of the following conditions holds:

(1) $r\leq3$;

(2) $r=\min\{d,n\}$;

(3) $r=4$ and $\min\{d,n\}\leq6$.
\end{thr}

\begin{proof}
Let $A\in\R^{d\times n}$. If $\trop(A)<\min\{d,n\}$, then Theorem~\ref{dssthr} implies that $\Kap(A)<\min\{d,n\}$.
If $\trop(A)<3$, then Theorem~\ref{dssthr} shows that $\Kap(A)=\trop(A)$.
If $\trop(A)=3$ and $\min\{d,n\}\leq6$, then from Theorems~\ref{cjrthr} and~\ref{6x6r3} it follows that $\trop(A)=\Kap(A)$.

Now assume that the conditions (1)--(3) fail to hold. It is enough to show how to construct a matrix
$B\in\R^{d\times n}$ such that $\trop(B)<r$, $\Kap(B)\geq r$.
Indeed, for $r=4$, $\min\{d,n\}\geq7$, we construct it via Example~\ref{dssex} and Theorem~\ref{extr}.
For $5\leq r\leq\min\{d,n\}-1$, $\min\{d,n\}\geq6$, we use Example~\ref{6x6} and apply Theorems~\ref{extr} and~\ref{extr2}.
\end{proof}

It is noted in Question~\ref{questi} that the $r$-by-$r$ minors of a $d$-by-$n$ matrix form a tropical
basis if and only if every $d$-by-$n$ matrix of tropical rank less than $r$ has the Kapranov
rank less than $r$. Thus Theorem~\ref{maintheoremofallthis} follows directly from Theorem~\ref{mainthr2}.

\section{Acknowledgements}
I owe my deepest gratitude to my family for their unflagging love, patience and invaluable support
throughout my life and during the course of this work in particular.
I would like to thank my scientific advisor Prof. Alexander E. Guterman for helpful discussions
and constant attention to my work.

\end{document}